\documentclass[11pt]{article}

\usepackage{amssymb,amsmath,latexsym}

\newtheorem{thm}{Theorem}[section]
\newtheorem{prop}[thm]{Proposition}
\newtheorem{remark}[thm]{Remark}
\newtheorem{defn}[thm]{Definition}
\newtheorem{cor}[thm]{Corollary}

\newtheorem{lemma}[thm]{Lemma}

\numberwithin{equation}{section}

\def\RR{{\bf R}}

\def\FF{{\cal F}}

\def\RR{{\bf R}}

\def\A{{\bf A}}
\def\S{{\bf S}}
\def\C{{\cal C}}
\def\NN{{\cal N}}
\def\wt{\widetilde}
\def\wh{\widehat}

\def\wt{\widetilde}
\def\pf{\noindent{\bf Proof.} }
\def\eps{{\varepsilon}}
\def\E{{\bf E}}
\def\P{{\bf P}}

\def\qed{{\hfill $\Box$ \bigskip}}

\begin{document}
\title{\bf 
\vspace{2truein}
Relative Fatou's Theorem for $(-\Delta)^{\alpha/2}$-harmonic Functions in Bounded $\kappa$-fat Open Sets
\thanks{This research
is supported in part by NSF Grant
DMS-0071486.
}
\vspace{.4truein}
}

\author{ {\bf Panki Kim}
\bigskip\\
Department of Mathematics\\
University of Illinois \\
Urbana, IL 61801, USA \bigskip\\
Email: pankikim@math.uiuc.edu \bigskip\\
Telephone number:  (217) 244-5930\\
Fax number: (217) 333-9576 
}
\maketitle
\newpage
\vspace*{1.5truein}
\begin{abstract}
Recently it was shown in Kim \cite{K} that Fatou's theorem for
transient censored $\alpha$-stable processes in a bounded $C^{1,1}$ open set is true.
Here we give a probabilistic proof of relative Fatou's theorem for
$(-\Delta)^{\alpha/2}$-harmonic functions (equivalently for symmetric $\alpha$-stable processes) in   bounded $\kappa$-fat open set where $\alpha \in (0,2)$. That is, if $u$ is positive $(-\Delta)^{\alpha/2}$-harmonic function in  a bounded $\kappa$-fat open set $D$ and $h$ is singular positive $(-\Delta)^{\alpha/2}$-harmonic function in $D$, then non-tangential limits of  $u/h$ exist almost everywhere with respect to the Martin-representing measure of $h$.
This extends the result of Bogdan and Dyda \cite{BD}.
It is also shown that, under the gaugeability assumption, relative Fatou's theorem is true for operators obtained from the generator of the killed $\alpha$-stable process in bounded $\kappa$-fat open set $D$ through non-local Feynman-Kac transforms. 
As an application, 
 relative Fatou's theorem for relativistic stable processes is also true if $D$ is bounded $C^{1,1}$-open set.
\end{abstract}
\vspace{1truein}

\noindent {\bf AMS 2000 Mathematics Subject Classification}: Primary
31B25, 60J75; Secondary 60J45, 60J50

\vspace{.3truein}

\noindent {\bf Keywords and phrases:}  Green function, symmetric
stable process, relative Fatou's theorem, Martin kernel, Martin boundary, 
harmonic function, Feynman-Kac transforms,
Martin representation

\vfill \eject

\begin{doublespace}

\section{Introduction}
Fatou \cite{F} in 1906 showed that bounded harmonic functions in the open unit disk have non-tangential limits almost everywhere on the unit circle. Later, Fatou's theorem (in classical sense) has been extended to some general open sets, up to uniform domains (see \cite{A2}, \cite{HW1}, \cite{HW2} and \cite{JK} for analytic approaches).
Probabilistic methods can be applied to proving Fatou's theorem. Through probabilistic methods, Fatou's theorem for Brownian motion (classical sense) and for various diffusion processes were proved (see, for example, \cite{B}, \cite{CDM}, \cite{D} and \cite{Du}). So far Fatou's theorem has mainly been established for elliptic differential operators or equivalently, diffusion processes. However, recently Fatou's theorem for discontinuous transient censored stable processes in bounded $C^{1,1}$ open set is proved in Kim \cite{K} through a probabilistic method.

Fatou's theorem can be stated in a more general setting, namely, relative Fatou's theorem. in 1959 Doob \cite{D3}  showed that 
the ratio $u/h$ of two positive harmonic functions for Brownian motion on an open solid sphere has non-tangential limits almost everywhere with respect to the Martin-representing measure of $h$ (see \cite{D2} for a non-probabilistic proof). Later, relative Fatou's theorem (in the classical sense) has been extended to some general open sets (for example, see \cite{W} and references therein). 

But relative Fatou's theorem stated above (and Fatou's theorem) is not true for $(-\Delta)^{\alpha/2}$-harmonic functions (see R. Bass and D. You \cite{BY} for some counterexamples). In this paper,  relative Fatou's theorem for $(-\Delta)^{\alpha/2}$-harmonic functions means the existence of non-tangential limits of  the ratio $u/h$ of positive $(-\Delta)^{\alpha/2}$-harmonic function $u$ in a open set $D$ and singular positive 
$(-\Delta)^{\alpha/2}$-harmonic function $h$ in $D$ (see Theorem \ref{T:Fatou} for the precise statement).

In \cite{BD}, K. Bogdan and B. Dyda prove relative Fatou's theorem for a special class of $(-\Delta)^{\alpha/2}$-harmonic functions in bounded $C^{1,1}$-domains through an analytic method. In this paper, through a probabilistic method, we show that relative Fatou's theorem for $(-\Delta)^{\alpha/2}$-harmonic functions as stated  in the previous paragraph holds for much more general open sets, namely, bounded $\kappa$-fat open sets, which include bounded Lipschitz open sets. The analogous result is unknown even in the Brownian motion case.

Since symmetric $\alpha$-stable process $X^D$ in an open set $D$ has discontinuous sample paths, there is
a large class of additive functionals of $X^D$ which are not
continuous.
Additive functionals
of the form 
\[
A_{q+F}(t)=\int_0^t q(X^D_s)ds + \sum_{s\le t}F(X^D_{s-}, X^D_s)
\]
constitute an important class of discontinuous additive functionals
of $X^D$.
Here  $q$ is a Borel measurable function on $D$ and
$F$ is some bounded Borel measurable function on $D \times D $ vanishing on
the diagonal. Such an additive functional defines a Feynman-Kac
semigroup by
$$ Q_t f(x) =\E_x \left[ \exp\left( A_{q+F} (t)\right) f(X^D_t) \right].
$$
 The above Feynman-Kac transform is called non-local (see Remark 1 of \cite{CS7}).
The Feynman-Kac transform of the above type has
been studied in \cite{C}, \cite{CK2}, \cite{CS6}, \cite{CS7} and \cite{CS4} in connection
with gauge and conditional gauge theorems for a large class
of Markov processes.
In this paper, we study the boundary behavior of the ratio of harmonic functions 
for a symmetric $\alpha$-stable process in a bounded $\kappa$-fat open set under possibly discontinuous Feynman-Kac perturbation.
Under the gaugeability assumption, 
we show that relative Fatou's theorem is also true under possibly discontinuous Feynman-Kac perturbation.
To our knowledge, the boundary behavior of harmonic functions
  under non-local Feynman-Kac
perturbations has not been studied previously except in Kim \cite{K}.

This paper is organized as follows. In section 2, we recall the definition of symmetric $\alpha$-stable process and collect some known facts concerning  symmetric $\alpha$-stable process in bounded $\kappa$-fat open set and $(-\Delta)^{\alpha/2}$-harmonic function from \cite{SW}.
Section 3 contains the proof of relative Fatou's theorem for $(-\Delta)^{\alpha/2}$-harmonic functions in  bounded $\kappa$-fat open sets. 
The main idea of our proof is similar to Kim \cite{K}, which is inspired by Doob's approach (see also Bass \cite{B}). We use Harnack and boundary Harnack principle obtained in \cite{SW} and extend some results in \cite{CS3} to bounded $\kappa$-fat open set. If the open set is the unit ball in $\RR^2$, we show that our result is the best possible.
In section 4, we recall the definition of new Kato classes from \cite{C}, and nonlocal Feynman-Kac transforms from \cite{C} and \cite{CK2}. Then under the gaugeability assumption, we show relative Fatou's theorem for non-local operators obtained from a  symmetric $\alpha$-stable process in bounded $\kappa$-fat open set through non-local Feynman-Kac transforms. As a consequence, relative Fatou's theorem for relativistic stable processes in bounded $C^{1,1}$ open set is established.

In this paper, we use ``$:=$" as a way of
definition, which is  read as ``is defined to be".
For functions $f$ and $g$, notation ``$f\approx g$"
means that there exist constants $c_2>c_1>0$ such that
$c_1 \, g \leq f \leq c_2 \, g$. The letter $c$, with or without subscripts, signifies a constant whose value is unimportant and which may change from location to location, even within a line.

The first version of this paper was finished in November of 2003. 
After this paper was submitted, we found that  
K. Michalik and M. Ryznar in \cite{MR} obtained 
the relative Fatou's theorem for singular $(-\Delta)^{\alpha/2}$-harmonic 
functions in bounded Lipschitz domains. As mentioned in \cite{MR}, 
the class of open sets considered in
the present paper is much broader than the class of bounded Lipschitz
domains (see below for the details). In fact, the open 
set considered in this paper does not need to be connected.
 The method in this paper is more probabilistic, while the
    \cite{MR} is more analytic in nature.
\medskip
\section{Preliminaries}

Let $X= \{X_t\}_{t \ge 0}$ denote
a symmetric $\alpha$-stable process
in $\RR^n$ with $\alpha \in (0, \, 2)$ and $n\geq 2$,
that is, let $X_t$ be
a L\'evy process
whose transition density
$p(t, y-x)$
relative to the Lebesgue measure is given
by the Fourier transform,
$$ \int_{\RR^n} e^{i x\cdot \xi}  p(t, x) dx
          =e^{-t |\xi|^\alpha} .
$$

Given  an open set $D\subset \RR^n$, define
$\tau_D:=\inf\{t>0: \, X_t\notin D\}$.
Let $X^D_t(\omega)=X_t(\omega)$ if $t< \tau_D(\omega)$
and set $X^D_t(\omega)=\partial$ if $t\geq  \tau_D(\omega)$,
where $\partial$ is a coffin state added to $\RR^n$. The process
$X^D$, i.e., the process $X$ killed upon leaving  $D$,
is called the (killed) symmetric $\alpha$-stable process in $D$.

To state Harnack principle for $X$,
we need the following definition.
\medskip
\begin{defn}\label{def:har1}
Let $D$ be an open subset of $\RR^n$.
 A locally integrable function
$u$ defined on $\RR^n$ taking values in $(-\infty, \, \infty]$
 and satisfying the condition 
$ \int_{\{x \in \RR^n ; |x| >1\}}|u(x)| |x|^{-(n+\alpha)} dx  <\infty $ is said to be

\begin{description}
\item{(1)} $(-\Delta)^{\alpha/2}$-harmonic in $D$ if
$$
\E_x\left[|u(X_{\tau_{B}})|\right] <\infty
\quad \hbox{ and } \quad 
u(x)= \E_x\left[u(X_{\tau_{B}})\right],
\qquad x\in B,
$$
for every open set $B$ whose closure is a compact
subset of $D$;

\item{(2)} $(-\Delta)^{\alpha/2}$-superharmonic in $D$ if
$u$ is lower semicontinuous in $D$ and
$$
\E_x\left[u^- (X_{\tau_{B}})\right] <\infty
\quad \hbox{ and } \quad 
u(x)\geq  \E_x\left[u(X_{\tau_{B}})\right], \qquad x\in B,
$$
for every open set $B$ whose closure is a compact
subset of $D$;

\item{(3})
 regular $(-\Delta)^{\alpha/2}$-harmonic in $D$ if it is $(-\Delta)^{\alpha/2}$-harmonic in $D$  and 
for each $x \in D$, 
$$
u(x)= \E_x\left[u(X_{\tau_{D}})\right];
$$

\item{(4})
 singular $(-\Delta)^{\alpha/2}$-harmonic in $D$ if it is $(-\Delta)^{\alpha/2}$-harmonic in $D$ and it vanishes outside $D$
\end{description}
\end{defn}
\medskip
Note that a $(-\Delta)^{\alpha/2}$-harmonic function
in an open subset $D$ is
continuous on $D$ (see \cite{BB} for an analytic definition and its equivalence). Also note that singular $(-\Delta)^{\alpha/2}$-harmonic function  $u$ in $D$
is harmonic with respect to $X^D$. i.e. 
$$
\E_x\left[|u(X^D_{\tau_{B}})|\right] <\infty
\quad \hbox{ and } \quad 
u(x)= \E_x\left[u(X^D_{\tau_{B}})\right],
\qquad x\in B,
$$ 
for every open set $B$ whose closure is a compact
subset of $D$.
\medskip

\begin{thm}\label{T:Har}
(Bogdan \cite{Bo})
Let $x_{1}, x_{2}\in D$, $r>0$ 
such that $|x_{1}-x_{2}|< Mr$. Then there exists a constant $J$ depending only
on $n$ and $\alpha$, such that
$$
J^{-1}M^{-(n+\alpha)}u(x_{2})\leq u(x_{1})
\leq JM^{n+\alpha}u(x_{2})\,
$$
for every nonnegative $(-\Delta)^{\alpha/2}$-harmonic function $u$ in
$B(x_{1}, r)\cup B(x_{2},r)$.
\end{thm}
\medskip
Unlike Brownian motion, the above theorem does not require Harnack chain argument. This observation is one of the reason why relative Fatou's theorem for $(-\Delta)^{\alpha/2}$-harmonic function is true in very general open set. 

We will often use the results in \cite{SW}. For our convenience, their main results are listed here. First we adopt the definition of $\kappa$-fat open set from \cite{SW}. 
\medskip  
\begin{defn}\label{fat}
Let $\kappa \in (0,1/2]$. We say that an open set $D$ in 
$\RR^n$ is $\kappa$-fat if there exist $R>0$ such that for each $z \in \partial D$ and $r \in (0, R)$, 
$D \cap B(z,r)$ contains a ball $B(a(r,z),\kappa r)$. The pair $(R, \kappa)$ is called the characteristics of 
the $\kappa$-fat open set $D$.
\end{defn}

\medskip

Note that every Lipschitz domain and
{\it non-tangentially accessible domain} defined by Jerison and Kenig in \cite{JK} are $\kappa$-fat. Moreover, 
every {\it John domain} is  $\kappa$-fat (see Lemma 6.3 in \cite{MV}).
The boundary of a $\kappa$-fat open set can be highly
nonrectifiable and, in general, no regularity of its boundary
can be inferred.  Bounded $\kappa$-fat open set can even be locally disconnected (see examples immediately following Lemma \ref{seq}). 

 Throughout this paper, $D$ is a bounded $\kappa$-fat open set, $n \ge 2$, $\alpha \in (0, 2)$ and ${\cal F}_t$ is the completed filtration for $X^D_t$, that is,
$$  \displaystyle{{\cal F}_t :=\bigcap_{x \in D} ~\sigma(\sigma(X^D_s : 0 \le s \le t) \cup \NN^x)} $$ where $\NN^x$ is the collection of $P_x$-null sets.
\medskip

\begin{thm}\label{T:BH}
(Theorem 3.1 in \cite{SW})
Let $D$ be a bounded open set in $\RR^n$ which is $\kappa$-fat for some $\kappa \in (0,1/2]$.
Then there exists constant $C=C(n, \alpha) >1$ 
such that for any $z \in \partial D$, $r \in (0,R)$ and functions $u,v \ge 0$ in $\RR^n$, regular 
  $(-\Delta)^{\alpha/2}$-harmonic in $D \cap B(z,2r)$, vanishing on $D^c \cap B(z,2r)$, we have
$$
C^{-1} \kappa^{n+ \alpha} \frac{u(a(r,z))}{v(a(r,z))} 
\le \frac{u(x)}{v(x)} \le C \kappa^{-n- \alpha}\frac{u(a(r,z))}{v(a(r,z))}, 
~~~~~ x\in D\cap B\left(z,\frac{r}{2}\right).
$$
\end{thm}

\medskip
It is well known that there is a positive continuous symmetric function
$G_D(x, y)$ on
$(D\times D)\setminus d$, where $d$ denotes the diagonal,
such that for any Borel measurable function $f\geq 0$,
$$ \E_x \left[ \int_0^{\tau_D} f(X_s) ds \right]
=\int_D G_D (x, y) f(y) \, dy.
$$
We set $G_D$ equal to zero on the diagonal of $D \times D$ and outside $D \times D$.
Function $G_D(x, y)$ is called the Green function of $X^D$,
or the Green function of $X$ in $D$.
For any $x \in D$, $G_D( \, \cdot \, ,x)$ is singular $(-\Delta)^{\alpha/2}$-harmonic in $D\setminus \{x\}$ 
and regular $(-\Delta)^{\alpha/2}$-harmonic in $D\setminus B(x, \eps)$ for every $\eps >0$.
 When $D=\RR^n$, it is well known that 
$$
G(x,y):=G_{\RR^n}(x,y) = A(n,\alpha)|x-y|^{\alpha-n}, ~~~~x,y \in \RR^n
$$
where $A(n,\alpha)$ is a positive constant depending only on $n$ and $\alpha$.

Fix $x_0 \in D$ and set
$$M_D(x,y) := \frac {G_D(x,y)}{G_D(x_0 , y)} ~,~~~~ x,y \in D.$$
It is shown in \cite{SW} that $M_D(x,z):= \lim_{y\rightarrow z \in \partial D} M_D(x,y)$ exists for every
$z\in \partial D$, which is called the Martin kernel
of $D$, and that $M_D(x, z)$ is jointly  continuous in
$D\times \partial D$. For each $z \in \partial D$, set $M_D(x,z)=0$ for $x \in D^c$.
The following properties of the Martin kernel and the Martin boundary   with respect to $X^D$ are established in \cite{SW}.
We call $w \in \partial D$ a regular boundary point if $\P_w(\tau_D =0 )=1$.
\medskip

\begin{thm}\label{T:M}(Theorem 4.1 in \cite{SW})
For each $z\in \partial D$,
$x\mapsto M_D(x, z)$ is a minimal singular $(-\Delta)^{\alpha/2}$-harmonic function, 
and
 the Martin boundary and the minimal Martin
boundary of $D$ can all be identified with the Euclidean boundary
$\partial D$ of $D$. Moreover, $M_D(x, z) \to 0$ as $D \ni x \to w$ for every regular boundary point $w \not= z$.  

\end{thm}

\medskip

 $X^D$ is a transient symmetric Hunt process satisfying Hypothesis (B)
in Kunita and Watanabe \cite{KW}. Thus non-negative singular $(-\Delta)^{\alpha/2}$-harmonic functions
admit a Martin representation.
Therefore, Theorem \ref{T:M} implies that, for every non-negative singular $(-\Delta)^{\alpha/2}$-harmonic function $u$,
there is a unique finite
measure $\nu$ on $\partial D$ such that
\begin{equation}\label{eqn:martin}
 u(x) = \int_{\partial D} M_D(x,z) \nu (dz),~~~~ x \in D.
\end{equation}
\medskip

\section{Relative Fatou's Theorem for $(-\Delta)^{\alpha/2}$-harmonic Functions
 }

In  this section, we establish relative Fatou's theorem for $(-\Delta)^{\alpha/2}$-harmonic functions
 in bounded $\kappa$-fat open set. 

The proof of the next proposition is well known (for example, see \cite{B} and \cite{K}).

\medskip
   
\begin{prop}\label{prop:3.2}
Given $ 0 < \lambda < 1$, there exists $c = c ( D, \alpha, \lambda)>1$ such that if $y \in D $ and $|y - x_0 | > 2 \delta_D (y) $ then 
\begin{equation}\label{eqn:3.5}
 \P_{x_0} \left(T_{B^{\lambda}_y}  < \tau_D \right) \ge c ~G_D(x_0, y) \delta_D (y)^{n- \alpha}
\end{equation}
where $B^{\lambda}_y := B(y, \lambda \delta_D (y))$, $ \delta_D (x)={\rm dist}(x,D^c)$ 
and $T_{B^{\lambda}_y} = \inf\{t>0: \, X_t \in B^{\lambda}_y \}$.
\end{prop}
\pf
First note that $x_0 \not\in B(y, \delta_D(y))$.
Since $G_D(x_0 , \,\cdot\,)$ is singular $(-\Delta)^{\alpha/2}$-harmonic in $D \setminus \{x_0\}$ , by Theorem \ref{T:Har}, there exists $c = c ( D, \alpha, \lambda)>1$ such that
\begin{equation}\label{eqn:3.6}
G_D1_{B^{\lambda}_y}(x_0) \ge c~ G_D(x_0 , y) \delta_D (y)^n .
\end{equation}
Using the strong Markov property, one can easily see that
\begin{equation}\label{eqn:3.7}
G_D1_{B^{\lambda}_y} (x_0) 
\le \P_{x_0}\left(T_{B^{\lambda}_y} <\tau_D \right) \sup_{w \in \overline{B^{\lambda}_y}} \E_w \int^{\tau_D}_0 1_{B^{\lambda}_y} (X_s) ds.
\end{equation}
Since 
\begin{equation}\label{eqn:3.8}
\E_w \int^{\tau_D}_0 1_{B^{\lambda}_y} (X_s) ds \le \int_{B^{\lambda}_y} G(w,v)dv \le c\int_{B^{\lambda}_y} \frac{dv}{|w-v|^{n- \alpha}}\le c \,\delta_D (y)^{\alpha}
\end{equation}
for every $w \in \overline{{B^{\lambda}_y}}$ where
$G(w,v)$ is Green function of $X$ in $\RR^n$,
we have (\ref{eqn:3.5}).
\qed

\medskip

Now let $(\P^z_x , X^z_t)$ be the $M_D(\cdot , z)$-transform of  $(\P_x , X^D_t)$, (killed) symmetric $\alpha$-stable process in $D$. That is, 
$$\P^z_x (A) := \E_x \left[ \frac{M_D(X^D_t , z)}{M_D(x,z)}; A \right]$$ 
if $A \in {\cal F}_t.$

First we show the following simple lemma, which is similar to Theorem 2.4 in \cite{CS3}.
\medskip

\begin{lemma}\label{lemma:boy}
For each $z \in \partial D$, 
 $M_D(\, \cdot \, , z)$ is bounded regular $(-\Delta)^{\alpha/2}$-harmonic 
in $D \setminus B(z, \eps)$ for every $\eps > 0$.
\end{lemma}

\pf
Fix $z \in \partial D$ and $\eps > 0$, and let $h(x):=M_D(x, z)$ for $x \in \RR^n$.
By Theorem \ref{T:BH},  $h$ is bounded on  $\RR^n \setminus B(z, \eps/2)$. In fact,
there exists $C=C(n, \alpha)>1$ such that for every $x \in D \setminus B(z, \eps/2)$,
\begin{eqnarray*}
&&M_D(x,z)\,=\,\lim_{D \ni y \to z} \frac{G_D(x,y)} {G_D(x_0,y)}\,
\le\, c \kappa^{-n- \alpha}\,  \frac{G_D(x,a)} {G_D(x_0,a)}\\
&&\le\, c \kappa^{-n- \alpha} \, \frac{G(x,a)} {G_D(x_0,a)}
\,\le\, C \kappa^{-n- \alpha}\,\sup_{y \in D \setminus B(z, \eps/2)}  \frac1{|y-a|^{n+\alpha}G_D(x_0,a)} < \infty
\end{eqnarray*}
 where $a=a(\eps/16,z)$ (see Definition \ref{fat}).
 
Take an increasing sequence of smooth open sets $\{ D_m \}_{m \ge 1}$ such that $\overline{D_m} \subset D_{m+1}$ and $\cup^{\infty}_{m=1} D_m = D \setminus B(z, \eps)$.
Set $\tau_m := \tau_{D_m}$ and $\tau_{\infty} := \tau_{D \setminus B(z, \eps)}$ .
Then $\tau_m \uparrow \tau_{\infty}$ and 
$\lim_{m \to \infty} X_{\tau_m} = X_{\tau_\infty} $ by quasi-left continuity of $X^D$.
Set $A = \{ ~\tau_m = \tau_\infty  ~\mbox{ for some } m \ge 1 \}$.
Let $N$ be the set of irregular boundary points of $D$. It is well known that Cap$(N)=0$ so that
\begin{equation}\label{zero1}
\P_x(X_{\tau_\infty} \in N)=0,  ~~~  x \in D.
\end{equation} 
We also know from \cite{SW} (Theorem \ref{T:M}) that if $w \in \partial D, w \not= z $
and $w$ is a regular boundary point, then 
$
h(x) \to 0 $ as  $ x \to w$
so that $h$ is continuous on $\overline{D \setminus B(z, \eps)} \setminus N$.
Therefore, since $h$ is bounded on $\RR^n \setminus B(z, \eps/2)$, 
by the bounded convergence theorem and (\ref{zero1}), we have
\begin{eqnarray*}
\lim_{m \to \infty} \E_x \left[ \,h (X^D_{\tau_m}) \,;\, \tau_m < \tau_\infty  \right] 
&=&\lim_{m \to \infty} \E_x \left[ \,h (X_{\tau_m})1_{ \overline{D \setminus B(z, \eps)} \setminus N}(X_{\tau_\infty}) \,;\, \tau_m < \tau_\infty \right]\\
&=& \E_x \left[ \,h (X_{\tau_\infty})1_{\overline{D \setminus B(z, \eps)} \setminus N}(X^D_{\tau_\infty}) \,;\, A^c\,\right]\\
&=& \E_x \left[ \,h (X_{\tau_\infty}) \,;\, A^c\,\right].
\end{eqnarray*}
By the boundedness of $h$ on  $\RR^n \setminus B(z, \eps/2)$, we can find two smooth open sets $U_1$ and $U_2$ such that $\overline{D \setminus B(z, \eps)} \subset U_1 \subset \overline{U}_1 \subset U_2$ and $h$ is the bounded on $\overline{U}_2$.
The rest of the proof is similar to the proof in Theorem 2.4 in \cite{CS3}. So we omit it here.
\qed

\medskip

The following theorem is proved in \cite{CS3} for bounded Lipschitz domain. In fact, by a similar proof with some modifications, it is true for bounded $\kappa$-fat open set $D$ too.
\medskip
\begin{thm}\label{thm:son1}
$$
\P^z_x \left( \lim_{t \uparrow \tau^z_D} X^z_t = z ,~ \tau^z_D< \infty \right) = 1 ~~~\mbox{ for every } x \in D,~ z \in \partial D.
$$
\end{thm}

\pf
By 3G theorem (Theorem 6.1 in \cite{SW}), there exists a positive constant $C=C(D,\alpha, n) < \infty$ such that
$$
\sup_{x \in D, z \in \partial D} \E^z_x[\tau^z_D] ~\le~ C. 
$$  
In fact, by Lemma 3.11 in \cite{CS3} and Theorem 6.1 in \cite{SW}, we have 
\begin{eqnarray*}
&&\E^z_x[\tau^z_D]\, = \,\E^z_x \int_{0}^{\infty} 1_{\{t < \tau^z_D\}}\,dt\,
\,=\,\frac1{M_D(x,z)} \int_{0}^{\infty} \E_x \left[M_D(X^D_t,z); \, t < \tau_D \right]dt \\
&&=\int_D\frac{G_D(x,y)M_D(y,z)}{M_D(x,z)} \,dy
\,\le\, 2\, c \,\sup_{x \in \overline{D}} \int_D |x-y|^{\alpha-n}\, dy ~< ~\infty. 
\end{eqnarray*}
Therefore 
$
\P^z_x ( \tau^z_D< \infty ) = 1$ for every  $x \in D$ and  $z \in \partial D.$

Now we fix $x \in D$ and $z \in \partial D$ and claim that $\P^z_x ( \lim_{t \uparrow \tau^z_D } X^z_t = z ) = 1$.
The proof of this claim is well-known (see Theorem 5.9 in \cite{CZ} and Theorem 3.17 in \cite{CS3}). 
We give a sketch of the proof here.

Let $r_m= 1/2^m$, $B_m:=B(z,r_m)$, $D_m:= D \setminus \overline{B_m}$ and set
$T_m := T_{B_m}$ and $R_m = \tau_{B_m \cap D}$.
We may suppose that $x \in D_m$.
By Lemma \ref{lemma:boy}, we have 
$$
M_D(x,z) \,=\,\E_x [ M_D ( X_{\tau_{D_m}}, z ) ] 
\,=\,\E_x [ M_D ( X_{T_m}, z ) ; T_m < \tau_D ] 
\,=\, M_D(x,z)\, \P^z_x (T^z_m < \tau^z_D) .
$$
It follows that for all $m \ge 1$ we have $\P^z_x (T^z_m < \tau^z_D) = 1$.
Let $L_k := \sup_{y \in B^c_k \cap D} M_D(y,z)$, which is finite by  Lemma \ref{lemma:boy}.
For $k < m$, 
$$
\P^z_x \left[ T^z_m < \tau^z_D ,~ R^z_k \circ \theta_{T^z_m} < \tau^z_D \right] ~\le~ \frac{L_k}{M_D(x,z)}\P_x(T_m < \tau_D)
$$
(see page 283 in \cite{CS3} for details).
Since $\{z\}$ in $\RR^n$ with $n \ge 2$ has zero capacity with respect to $X$ (for example, see \cite{BBC} for details), we have
$$
\limsup_{m \to \infty} \P_x (T_m < \tau_D) ~ \le~\P_x (T_{\{z\}} \le \tau_D ) ~\le~ \P_x (T_{\{z\}} < \infty ) ~=~ 0.
$$
The rest of the proof is the same as the proof in Theorem 3.17 in \cite{CS3}.
\qed
\medskip

The above Theorem implies  
that 
$\P^{\,\cdot}_x ( \lim_{t \uparrow \tau_D} X_t \in K ) = 1_K(\,\cdot\,)$ ~for every $x \in D$ and Borel subset $K \subset \partial D$.
So the next theorem, which is stated for bounded Lipschitz domain in \cite{CS3}, follows easily. 
For positive singular $(-\Delta)^{\alpha/2}$-harmonic function $h$ in $D$, we let $(\P^h_x , X^h_t)$ be the $h$-transform of  $(\P_x , X^D_t)$, that is, 
$$\P^h_x (A) := \E_x \left[ \frac{h(X^D_t)}{h(x)}; A \right]$$ 
if $A \in {\cal F}_t.$

\medskip

\begin{thm}\label{thm:son2}
Let $\nu$ be a finite measure on $\partial D$.
Define
$$
h(x) := \int_{\partial D} M_D(x, w) ~\nu (dw), ~~~~ x \in D,
$$ 
which is a positive singular $(-\Delta)^{\alpha/2}$-harmonic function in $D$. 
Then for $ x \in D$,
$$
\P^h_x \left(\lim_{t \uparrow \tau^h_D} X^h_t \in K\right) = \frac1{h(x)} \int_K 
M_D(x,w) \nu(dw)
$$
 where $K$ is a Borel subset of $\partial D$.
\end{thm}
\medskip
The  following result is a easy consequence of the above theorem.
\medskip

\begin{prop}\label{prop:3.3}
Let $\nu$ be a finite measure on $\partial D$ with $\nu(\partial D)=1$.
Define
$$
h(x) := \int_{\partial D} M_D(x, w) ~\nu (dw), ~~~~x \in D,
$$
 so that $h(x_0)=1$. If $A \in {\cal F}_{\tau_D}$, then  
$$
\int_K \P^z_{x_0} (A)\nu(dz)
= \P^h_{x_0}\left(A \cap \left\{\lim_{t \uparrow \tau^h_D} X^h_t \in K\right\} \right).
$$
for every Borel subset $K$ of $\partial D$.
\end{prop}
\pf
Take an increasing sequence of open sets $\{D_m\}_{m \geq 1}$ such that $\overline{D_m} \subset D_{m+1}$ and $\cup^{\infty}_{m=1} D_m = D$.
Set $\tau_m = \tau_{D_m}$ and fix an $A \in {\cal F}_{\tau_{m}}$. Since $M(x_0 , z)=1$ for $z \in \partial D$, by Theorem \ref{thm:son2}, Fubini's Theorem and the strong Markov property for conditional process (for example, see \cite{KW}), we have that for every Borel subset $K$ of $\partial D$, 
\begin{eqnarray*}
\int_K \P^z_{x_0} (A)\nu(dz)  
&=& \int_K \E_{x_0} \left[ M_D(X^D_{\tau_m} , z )\, ;\, A\right]  \nu(dz)  \\
&=& \E_{x_0} \left[ \int_K M_D(X^D_{\tau_m}, z ) \nu( dz)  ;~ A \right] \\
&=& \E_{x_0} \left[h\left(X^D_{\tau_m}\right)\P^h_{X_{\tau_m}}\left(\lim_{t \uparrow \tau^h_D} X^h_t \in K\right) \,;\, A \right]\\
&=& \E^h_{x_0} \left[\P^h_{X^h_{\tau_m}}\left(\lim_{t \uparrow \tau^h_D} X^h_t \in K\right) ; A \right]\\
&=& \P^h_{x_0}\left(A \cap \left\{\lim_{t \uparrow \tau^h_D} X^h_t \in K \right\}\right).
\end{eqnarray*}
Let $m \rightarrow \infty$. Then
\begin{equation}\label{eqn:AA}
\int_K \P^z_{x_0} (A)\nu(dz)  
= \P^h_{x_0}\left(A \cap \left\{\lim_{t \uparrow \tau^h_D} X^h_t \in K\right\}\right)
\end{equation}
for every Borel subset $K$ of $\partial D$ and $A \in \cup_{m \geq 1}{\cal F}_{\tau_{m}}$.
By monotone class theorem, (\ref{eqn:AA}) is true for every Borel subset $K$ of $\partial D$ 
and $A \in {\cal F}_{\tau_D}$.
\qed

\begin{defn}\label{def:3.4}
$A \in {\cal F}_{\tau_D}$ is shift-invariant if whenever $T < \tau_D $ is a stopping time, $ 1_A \circ \theta_T = 1_A$ $\P_x$-a.s. for every $x \in D$.
\end{defn}

\medskip

The next proposition is well-known (see page 196 in Bass \cite{B}).

\medskip
\begin{prop}\label{prop:3.5}
(0-1 law) If $A$ is shift-invariant, then $x \to \P^z_x (A)$ is a constant function which is either $0$ or $1$.
\end{prop}

\pf
For every stopping time $T < \tau_D$,
$$
\P^z_x (A) \,=\,\P^z_x (A \circ \theta_T ) \,=\, \E^z_x \,\P^z_{X^z_T} (A)\, 
=\, \frac1{M_D(x,z)}\, \E_x \left[M_D(X^D_T , z) \P^z_{X^D_T} (A) \right] .
$$
 So $ M_D(\,\cdot\,, z) \P^z_{\cdot} (A)$ is positive singular $(-\Delta)^{\alpha/2}$-harmonic and bounded above by $M_D(\,\cdot \,, z )$.
Therefore $\P^z_{\cdot} (A)$ is constant by Theorem \ref{T:M}. With $D_m$ and $\tau_m$ in the proof of Proposition \ref{prop:3.3} and $B \in \FF_{\tau_m}$,
$$
\P^z_x (A \cap B) \,=\, \E^z_x \left[\P^z_{X^z_{\tau_m}} (A)\,;\, B \right] \,=\, \P^z_x (A)\, \P^z_x (B).$$
Let $m \rightarrow \infty$ and let $B =A$. Then we have $\P^z_x (A) = (\P^z_x (A))^2$. Therefore, $\P^z_x (A) = 0$ or $1$ for every $x \in D$.  \qed

\medskip

Now we define the Stolz open set for $\kappa$-fat open set $D$. Recall the characteristics $(R,\kappa)$ of $D$ from Definition \ref{fat}.
\begin{defn}\label{Stolz} 
For $z \in \partial D$ and $\beta > (1-\kappa)/\kappa$, let 
$$
A^{\beta}_z := \left\{ y \in D ;~ \delta_D(y) < \min\left\{\frac{\delta_D (x_0)}{3}, R\right\} ~\mbox{ and } ~|y-z| <  \beta \delta_D(y) \right\}.
$$
We call $A^{\beta}_z$ the Stolz open set for $D$ at $z$ with the angle $\beta > (1-\kappa)/\kappa$.
\end{defn}
\medskip
Since $\beta > (1-\kappa)/\kappa$,  $A^{\beta}_z$ is not empty. In fact, the following is true.
\medskip
\begin{lemma}\label{seq}
For every $z \in \partial D$ and $\beta > (1-\kappa)/\kappa$, there exists a sequence $\{y_k\}_{k \ge 1} \subset A^{\beta}_z$ 
such that $\lim_{k \to \infty} y_k =z$.
\end{lemma}
\pf
Fix $z \in \partial D$ and $\beta > (1-\kappa)/\kappa$. For $k \ge 1$, let 
$$
y_k := a\left(\frac{R}{2^k},z \right) 
~~\mbox{ so that } B\left(y_k, \frac{\kappa R}{2^k}\right) \,\subset\, B \left(z, \frac{R}{2^k}\right) \cap D.
$$
We may assume 
$$
\delta_D(y_k) < \min\left\{\frac{\delta_D (x_0)}{3}, R\right\}  ~~~~\mbox{ for every } k \ge 1.
$$
Since 
$$
\frac{R}{2^k} \,\ge\, |y_k-z| + \frac{\kappa R}{2^k}~~ \mbox{ and }~~ 
\delta_D(y_k) \,\ge\, \frac{\kappa R}{2^k}   ,
$$
we have
$$
|y_k-z|~ \le~ \frac{1- \kappa}{\kappa}\, \delta_D(y_k) ~ <~ \beta \delta_D(y_k). 
$$
Therefore  $\{y_k\}_{k \ge 1} \subset A^{\beta}_z$ . Clearly   $\lim_{k \to \infty} y_k =z$ because  $|y_k-z| \le R/2^k$.
\qed

\medskip

A simple example of bounded $\kappa$-fat open set is 
$$
G :=\left\{(x,y) \in \RR^2 \,;~ -1 < x < 1 , ~ 0 < y < 1 \right\}~ \setminus ~\bigcup_{k=1}^{\infty} E_k
$$
where 
$E_k:=\{ y = 2^{-k} \} $. One can easily see that every $B(z, \eps) \cap G$  is not locally connected for every $z \in \overline{G} \cap \{ y=0 \}$ and $\eps >0$.
In particular,  $A^{\beta}_z$ is  not locally connected at every $z \in \overline{G} \cap \{ y=0 \}$.
\medskip 
\begin{prop}\label{prop:3.6}
Given $z \in \partial D$ and $\beta > (1-\kappa)/\kappa$, there exists $c = c(D, \alpha , \lambda, x_0, \beta) > 0$ such that if $y \in A^{\beta}_z$ 
then 
$$
\P^z_{x_0} \left(T^z_{B^{\lambda}_y} < \tau^z_D \right) ~>~ c  
$$
where $B^{\lambda}_y := B(y, \lambda \delta_D (y))$, $ \delta_D (x):={\rm dist}(x,D^c)$ 
and $T^z_{B^{\lambda}_y} := \inf\{t>0: \, X^z_t \in B^{\lambda}_y \}$.
\end{prop}
\pf
Fix $z \in \partial D$ and  $\beta > (1-\kappa)/\kappa$.
Since $M_D(\cdot,z)$ is $(-\Delta)^{\alpha/2}$-harmonic function in $D$, by the Harnack principle (Theorem \ref{T:Har}) and Proposition \ref{prop:3.2} we have 
\begin{eqnarray*}
\P^z_{x_0} \left( T^z_{B^{\lambda}_y} < \tau^z_D \right) 
&=& \E_{x_0} \left[ M_D(X_{T_{B^{\lambda}_y}} , z) \,;\,  T_{B^{\lambda}_y} < \tau_D \right]\\ 
&\geq& c~ \P_{x_0} \left(  T_{B^{\lambda}_y} < \tau_D \right)~ M_D(y,z) \\
&\geq& c \, G_D(x_0, y) ~\delta_D (y)^{n-\alpha}~ M_D(y,z)\\
&=& c \, G_D(x_0, y) ~\delta_D (y)^{n-\alpha}~ \lim_{w \in D \rightarrow z} \frac{G_D(y,w)}{G_D(x_0,w)}.
\end{eqnarray*}
Since $\min\{|y-z|,|x_0-z|\} > \delta_D (y)/2$,
 by boundary Harnack principle (Theorem \ref{T:BH}),
$$
\frac{G_D(y,w_1)}{G_D(x_0,w_1)} \approx \frac{G_D(y,w_2)}{G_D(x_0,w_2)}
~~~~\mbox{ for every } w_1, w_2 \in B\left(z, \frac{\delta_D (y)}{8}\right)\cap D.
$$
Let $z_1:=a(\delta_D (y)/8,z)$ so that  
$$
B\left(z_1, \frac{\kappa \delta_D (y)}{8}\right) 
\,\subset\, 
B \left(z, \frac{ \delta_D (y)}{8}\right) \cap D
$$
 and  fix a point $z_2$ in $ \partial B(y, \delta_D (y)/8)$.
We see that
\begin{eqnarray*}
&& \delta_D(z_1) \,\ge\, \frac{\kappa \delta_D (y)}{8} \,>\,\frac{\delta_D (y)}{8(\beta+1)}\,,
~~\delta_D(z_2) \,>\,  \frac{\delta_D (y)}{2}\,,~~
|z_2- y| \,=\,\frac{\delta_D (y)}{8},\\
&&|z_2-x_0| \,\ge\, |x_0-y|-|y-z_2| \,\ge\, \delta_D (x_0) -\delta_D (y) - \frac{\delta_D (y)}{16}\,>\, \delta_D (y),\\
&&|z_1-x_0| \,\ge\, |x_0-z|-|z-z_1| \,\ge\, \delta_D (x_0) -\frac{\delta_D (y)}{16}
 \,>\, \delta_D (y),\\
&\mbox{ and }&|z_1-y| \,\ge\,  |y-z|-|z_1-z| \,\ge\, \delta_D (y)- \frac{\delta_D (y)} {8} \,>\,  \frac{\delta_D (y)}{2}. 
\end{eqnarray*}
Thus $G_D(y,\cdot)$ and
$G_D(x_0,\cdot)$ are  $(-\Delta)^{\alpha/2}$-harmonic functions  in  
$$
B\left(z_1,\frac{\delta_D (y)}{8(\beta+1)}\right) 
\cup 
B\left(z_2,\frac{\delta_D (y)}{8(\beta+1)}\right).
$$
Since
$$
|z_1-z_2| \,\le \,|z_1-z|+|z-y|+|y-z_2|
\,< \,\frac{\delta_D (y)}{8}+\beta \delta_D (y)+\frac{\delta_D (y)}{8},
$$
by Theorem \ref{T:Har}, we have 
$G_D(y,z_1)\approx G_D(y,z_2)$ and
$G_D(x_0,z_1)\approx G_D(x_0,z_2)$ .
Therefore
$$
\frac{G_D(y,w)}{G_D(x_0,w)} \approx \frac{G_D(y,z_2)}{G_D(x_0,z_2)}
~~~~\mbox{ for every } w \in B\left(z, \frac{\delta_D (y)}{8}\right)\cap D.
$$
On the other hand, by  Theorem \ref{T:Har} we have 
$G_D(x_0,z_2)\approx G_D(x_0,y)$.
Now by Green function estimate for ball and monotonicity of Green function for $X$, 
$$
G_D(y, z_2) \,\ge\, G_{B(z_2, \delta_D (y)/2)}(y, z_2) \,\ge\, c |y- z_2|^{\alpha -n} \,\ge\, c \delta_D (y)^{\alpha -n}.
$$
Thus
$$
G_D(x_0, y) \, \delta_D (y)^{n-\alpha}\, \lim_{D \ni w \to z} \frac{G_D(y,w)}{G_D(x_0,w)}~ >~c ~>~0.
$$
\qed

\medskip

The next proposition is a variation of Theorem \ref{T:Har}.

\medskip

\begin{prop}\label{prop:3.7}
Given $\eps >0$, there exists $ \lambda_0 = \lambda_0 (\eps, \alpha, n) \in (0,1)$ such that whenever $u$ is a positive $(-\Delta)^{\alpha/2}$-harmonic function in $D$,
$$
\frac{1}{1+ \eps }u(z)~ \le~ u(y)~ \le~ (1+ \eps ) u(z)
$$
 for every $y \in D$ and $ z \in B^{\lambda_0}_y:= B(y, \lambda_0 \delta_D (y)).$
\end{prop}
\pf
Fix $y \in D$ and let $\tau_{\lambda} := \tau_{B^{\lambda}_y}$. 
  By estimates of Poisson kernels of $X_t$ and the positivity of $u$
  (cf. \cite{CS1}), for every $z, w \in B^{\lambda_0}_y$ where  $ 0 < \lambda_0 < \lambda <1$,
$$
\E_w [ u(X_{\tau_{\lambda}} )] ~\le ~{\left(\frac{\lambda^2}{\lambda^2-{\lambda_0}^2}\right)}^{\alpha/2} {\left(\frac{\lambda+\lambda_0}{\lambda-\lambda_0}\right)}^{n}\E_z [ u(X_{\tau_{\lambda}} )].$$
Since $u$ is $(-\Delta)^{\alpha/2}$-harmonic in $D$, for every $ z,w \in B^{\lambda_1}_y$ we have
$$
u(w)
~ =~ \E_w [ u(X_{\tau_{\lambda}} )] 
~\le ~{\left(\frac{\lambda^2}{\lambda^2-{\lambda_0}^2}\right)}^{\alpha/2}  {\left(\frac{\lambda+\lambda_0}{\lambda-\lambda_0}\right)}^{n} u(z)
 ~\le~ (1 + \eps ) u(z) 
$$
if $\lambda_0 > 0$ is small. 
In particular, 
$$
\frac{1}{1+ \eps } u(z)~\le~ u(y) ~\le ~(1+ \eps ) u(z)
$$
 for every $ z \in B^{\lambda_0}_y$.
\qed

\medskip

Before proving relative Fatou's theorem for $(-\Delta)^{\alpha/2}$-harmonic function in $D$, we need the following definition.
\bigskip

\begin{defn}\label{D:5.1}
A nonnegative Borel measurable function
$f$ defined on $D$
is said to be

\begin{description}
\item{(1)} excessive with respect to $X^D$ if for every $x\in D$ and $t>0$,
$$
\E_x\left[  f(X^D_t)\right] \,\le\, f(x)~~~~ \mbox{and}~~~~
\lim_{t \downarrow 0} \E_x\left[  f(X^D_t)\right] \,=\, f(x).
$$

\item{(2)} superharmonic with respect to $X^D$ if
$f$ is lower semi-continuous in $D$,
and $$
f(x)\,\ge \,\E_x\left[  f(X^D_{\tau_{B}})\right],
\qquad x\in B,
$$
for every open set $B$ whose closure is a compact subset
of $D$.
\end{description}
\end{defn}
It is well known that $f$ is excessive with respect to $X^D$ if and only if $f$ is superharmonic with respect to $X^D$ (in fact, this result is true in much more general setting, see \cite{CK2}).

Now we are ready to show relative Fatou's theorem for $(-\Delta)^{\alpha/2}$-harmonic function in $D$. The proof is similar to the proof of Theorem 3.10 in \cite{K} but we spell out detail for the reader's convenience.

\medskip

\begin{thm}\label{T:Fatou}
Let $h$ be a positive singular $(-\Delta)^{\alpha/2}$-harmonic function in $D$ with the Martin-representing measure $\nu$.
That is, 
$$
h(x) = \int_{\partial D} M_D(x, w) \,\nu (dw),~~~ x \in D
$$ 
where $\nu$ is a finite measure on $\partial D$.
If $u$ is a nonnegative $(-\Delta)^{\alpha/2}$-harmonic function in $D$, then for $\nu$-a.e. $z \in \partial D$,
\begin{equation}\label{eqn:3.9}
\lim_{ A^{\beta}_z \ni x \rightarrow z} \frac{u(x)}{h(x)} \mbox{ exists for every } \beta >\frac{1-\kappa}{\kappa}.
\end{equation}
\end{thm}

\pf
Without loss of generality, we assume  $\nu(\partial D) = 1$.
Since $u$ is  non-negative superharmonic with respect to $X^D$,  $u$ is excessive with respect to $X^D$.
In particular, $\E_x[  u(X^D_t)] \le u(x)$ for every x $\in D$. So by Markov property for conditional process (for example, see \cite{SW}), we have for every $t, s >0$
$$
\E^h_{x_0}\left[  \frac{u(X^h_{t+s})}{h(X^h_{t+s})} \,\big|\, {\cal F}_s\right]
~=~ \E^h_{X^h_s}\left[  \frac{u(X^h_{t})}{h(X^h_{t})}\right]
~=~\frac1{h(X^h_s)}\E_{X^h_s}\left[  u(X^D_{t})\right]~\le~ \frac{u(X^h_s)}{h(X^h_s)}
$$
Therefore, we see that $u(X^h_t)/h(X^h_t)$ is a non-negative supermartingale with respect to $\P^h_{x_0}$. therefore the martingale convergence theorem gives
$$\lim_{t \uparrow \tau^h_D} \frac{u(X^h_t)}{h(X^h_t)} \mbox { exists and is finite }\P^h_{x_0}\mbox{-a.s. }, $$
so by Proposition \ref{prop:3.3}, we have 
\begin{eqnarray*}
1 &=& \P^h_{x_0} \left( \lim_{t \uparrow \tau^h_D} \frac{u(X^h_t)}{h(X^h_t)} \mbox { exists and is finite} \right) \\
&=& \int_{\partial D} \P^z_{x_0} \left(\lim_{t \uparrow \tau^z_D} \frac{u(X^z_t)}{h(X^z_t)} \mbox { exists and is finite}\right)~ \nu(dz).
\end{eqnarray*}
Therefore, for $\nu$-a.e. $z \in \partial D$ 
\begin{equation}\label{eqn:3.10}
\P^z_{x_0} \left(\lim_{t \uparrow \tau^z_D} \frac{u(X^z_t)}{h(X^z_t)} \mbox { exists and is finite} \right) = 1.
\end{equation}

We are going to show that (\ref{eqn:3.9}) holds for $z \in \partial D$ satisfying (\ref{eqn:3.10}).
Fix $z \in \partial D$ satisfying (\ref{eqn:3.10}) and fix a $\beta >(1-\kappa)/\kappa$.
Let 
$$
l ~:=~ \limsup_{ A^{\beta}_z \ni y \rightarrow z} \frac{u(y)}{h(y)} 
$$
 and assume $l < \infty.$
Then by Lemma \ref{seq}, for any $\eps > 0$ there exists a sequence $\{ y_k \}^{\infty}_{k=1} \subset A_z^{\beta} $ such that $u(y_k)/h(y_k) > l/(1+\eps )$ and $y_k \rightarrow z.$
By Proposition \ref{prop:3.7}, there is $\lambda_0 = \lambda_0 (\eps, \alpha, n)>0 $ such that 
\begin{equation}\label{eqn:3.11}
\frac{u(w)}{h(w)} ~\geq~ \frac{u(y_k)}{(1+\eps )^2h(y_k)}~ >~ \frac{l}{(1+\eps )^3}
\end{equation}
for every $ w \in B^{\lambda_0}_{y_k}=B(y_k, \lambda_0 \delta_D (y))$.

On the other hand, 
$$
\P^z_{x_0} \left(T^z_{B^{\lambda_0}_{y_k}} < \tau^z_D ~\mbox { i.o.} \right)
~ \geq~ \liminf_{k \rightarrow \infty} \,\P^z_{x_0} \left(T^z_{B^{\lambda_0}_{y_k}} < \tau^z_D \right) ~\geq~ c~  >~ 0 .
$$
But $\{ T^z_{B^{\lambda_0}_{y_k}} < \tau^z_D ~\mbox { i.o.} \}$ is shift-invariant.
Therefore by Proposition \ref{prop:3.5}
\begin{equation}\label{eqn:3.12}
\P^z_{x_0} \left( X^z_t \mbox{ hits infinitely many } B^{\lambda_0}_{y_k} \right) 
~=~\P^z_{x_0} \left(T^z_{B^{\lambda_0}_{y_k}} < \tau^z_D ~~\mbox { i.o.} \right) 
~=~1 .
\end{equation}
From (\ref{eqn:3.10})-(\ref{eqn:3.12}), we have 
$$
\lim_{t \uparrow \tau^z_D} \frac{u(X^z_t)}{h(X^z_t)} 
~\geq~ \frac{l}{(1+\eps )^3} ~~~~\P^z_{x_0} \mbox{-a.s. for every } \eps > 0 .$$
Letting $\eps \downarrow 0$, 
\begin{equation}\label{eqn:3.13}
\lim_{t \uparrow \tau^z_D} \frac{u(X^z_t)}{h(X^z_t)} 
~\geq ~\limsup_{ A^{\beta}_z \ni y \rightarrow z} \frac{u(y)}{h(y)} ~~~~\P^z_{x_0} \mbox{-a.s. }.
\end{equation}
If $~l=\infty$, then for any $M > 1$, there exists a sequence $\{ y_k \}^{\infty}_{k=1} \subset A_z^{\beta} $ such that $u(y_k)/h(y_k)> 4 M$ and $y_k \rightarrow z.$
By Proposition \ref{prop:3.7}, there is $\lambda_1 = \lambda_1 (M, \alpha, n) >0 $ such that 
$$
\frac{u(w)}{h(w)} ~\geq~ \frac{M^2u(y_k)}{(M+1)^2 h(y_k)} ~> ~M
$$
for every $ w \in B^{\lambda_1}_{y_k}$.
So similarly we have
$$
\lim_{t \uparrow \tau^z_D} \frac{u(X^z_t)}{h(X^z_t)}~ >~ M  ~~~~\P^z_{x_0} \mbox{-a.s. }
$$
for every $M >1$, which is a contradiction because $\lim_{t \uparrow \tau^z_D} u(X^z_t)/h(X^z_t)$ is finite $\P^z_{x_0}$-a.s.. Therefore  $l < \infty$.

Now let 
$$m := \liminf_{ A^{\beta}_z \ni y \rightarrow z} \frac{u(y)}{h(y)} < \infty.$$
Then for any $\eps > 0$, there exists a sequence $\{ z_k \}^{\infty}_{k=1} \subset A_z^{\beta} $ such that $u(z_k)/h(z_k)< m(1+\eps )$ and $z_k \rightarrow z.$
By Proposition \ref{prop:3.7},  
\begin{equation}\label{eqn:3.11A}
\frac{u(w)}{h(w)} ~\leq~ (1+\eps )^2\,\frac{u(z_k)} {h(z_k)} ~<~ (1+\eps )^3m
\end{equation}
for every $ w \in B^{\lambda_0}_{z_k}$.
Similarly we have
\begin{equation}\label{eqn:3.12A}
\P^z_{x_0} \left( X^z_t \mbox{ hits infinitely many } B^{\lambda_0}_{z_k} \right) 
=1 .
\end{equation}
From (\ref{eqn:3.10}), (\ref{eqn:3.11A}) and (\ref{eqn:3.12A}), by letting $\eps \downarrow 0$ we have 
\begin{equation}\label{eqn:3.14}
\lim_{t \uparrow \tau^z_D} \frac{u(X^z_t)}{h(X^z_t)} 
~\leq~ \liminf_{ A^{\beta}_z \ni y \rightarrow z} \frac{u(y)}{h(y)} ~~~~\P^z_{x_0} \mbox{-a.s. }.
\end{equation}

We conclude from (\ref{eqn:3.13}) and (\ref{eqn:3.14}) that 
$$\lim_{ A^{\beta}_z \ni y \rightarrow z} \frac{u(y)}{h(y)} \mbox{ exists and is finite for }\nu\mbox{-a.e. } z \in \partial D.$$
\qed 

\medskip

\begin{remark}
{\rm Since constant functions in $\RR^n$ are $(-\Delta)^{\alpha/2}$-harmonic in $D$, one can easily see that the above Theorem is also true for every $(-\Delta)^{\alpha/2}$-harmonic function $u$ in $D$ either bounded from below or above.
}
\end{remark}
\medskip

\begin{remark}
{\rm If $D$ is a bounded $\kappa$-fat open set with $\sigma(\partial D) < \infty$ where $\sigma$ is the surface measure on $\partial D$, then for every  $(-\Delta)^{\alpha/2}$-harmonic function $u$ in $D$ either bounded from below or above,
\begin{equation}\label{FF}
\lim_{ A^{\beta}_z \ni x \rightarrow z} ~\frac{u(x)}{\int_{\partial D} M_D(x, w) \sigma (dw)}~ \mbox{ exists and is finite for }\sigma\mbox{-a.e. } z \in \partial D.
\end{equation}
In particular, (\ref{FF}) is true for bounded Lipschitz open sets. 
Therefore Theorem \ref{T:Fatou} extends the result of Bogdan and Dyda \cite{BD}.
}
\end{remark}

\medskip
With an extra condition, we can state relative Fatou's theorem for (non-singular) $(-\Delta)^{\alpha/2}$-harmonic functions.
\medskip
\begin{cor}\label{C1:Fatou}
Let $\nu$ be a finite measure on $\partial D$ and let
$$
h(x) := \int_{\partial D} M_D(x, w) \,\nu (dw), ~~~~ x \in D
$$
 If $v$ is a positive $(-\Delta)^{\alpha/2}$-harmonic function in $D$ and 
for $\nu$-a.e. $z \in \partial D$
\begin{equation}\label{TT}
\lim_{ A^{\beta}_z \ni x \rightarrow z} \frac{v(x)}{h(x)} \not= 0~~~~ \mbox{ for } \nu\mbox{-a.e. } z \in \partial D
 \mbox{ and every } \beta>\frac{1-\kappa}{\kappa}
\end{equation}
then for every positive $(-\Delta)^{\alpha/2}$-harmonic function $u$ in $D$,
$$
\lim_{ A^{\beta}_z \ni x \rightarrow z} \frac{u(x)}{v(x)}~ \mbox{ exists for } \nu\mbox{-a.e. } z \in \partial D \mbox{ and every } \beta
 >\frac{1-\kappa}{\kappa}.
$$
\end{cor}
\medskip
In particular, if  $D$ is a bounded $C^{1,1}$-open set and the Martin-representing measure of $h$ is  the surface measure on $\partial D$, 
we can replace the condition (\ref{TT}) to a concrete one.
\begin{cor}\label{C2:Fatou}
Let $D$ be a bounded $C^{1,1}$-open set
and $\sigma$ be the surface measure on $\partial D$.
 If $v$ is a positive $(-\Delta)^{\alpha/2}$-harmonic function in $D$ and 
there exists $c >0$ such that
$v(x) \ge c \delta_D(x)^{\alpha/2-1} $ for $x \in D$, then for every positive $(-\Delta)^{\alpha/2}$-harmonic function $u$ in $D$,
$$
\lim_{ A^{\beta}_z \ni x \rightarrow z} \frac{u(x)}{v(x)} \mbox{ exists for } \sigma \mbox{-a.e. } z \in \partial D \mbox{ and every } \beta > 1.
$$
\end{cor}
\pf
Let
$$
h(x) := \int_{\partial D} M_D(x, w) \,\sigma (dw), ~~~~ x \in D
$$
where $\sigma$ is the surface measure on $\partial D$.
By Corollary \ref{C1:Fatou}, It is enough to show that there exists $c >0$ such that
$h(x) \le c \delta_D(x)^{\alpha/2-1} $ for $x \in D$.
Since $D$ is bounded $C^{1,1}$ open set, there exists $c_1$ depending only on $D$ such that for every $x \in D$ and $k \ge 1$,
\begin{equation}\label{aA1}
\sigma ( \partial D \cap B (z_x, 2^k \delta_D (x)) )~ \le~ c_1(2^k \delta_D (x) )^{n-1}
\end{equation}
where $z_x \in \partial D$ and $|z_x - x| = \delta_D (x)$.
Let 
\begin{eqnarray*}
&&E_1 := \left\{ w \in \partial D\, ;~ |w - z_x| \le 2 \delta_D (x) \right\} \\
\mbox{and}
&&E_k := \left\{ w \in \partial D\, ;~ 2^{k-1} \delta_D (x) < |w - z_x| \le 2^k \delta_D (x) \right\} ~~\mbox{ for } k \ge 2.
\end{eqnarray*}
We see that
 for $k \ge 1$
\begin{equation}\label{aA3}
|x - w| ~\ge~ 2^{k -2} \delta_D(x) ~~~\mbox{ if } w \in E_k  .
\end{equation}
By (\ref{aA1}), (\ref{aA3}) and the Martin kernel estimate (Theorem 3.9 in \cite{CS3}) , 
\begin{eqnarray*}
h(x)&\le& c ~\delta_D(x)^{\alpha/2} \int_{\partial D} \frac
{\sigma (dw)}{|x-w|^{n}} \\
&=& c ~\delta_D(x)^{\alpha/2} \sum^{\infty}_{k=1} \int_{E_k} \frac{\sigma (dw)}{|x-w|^{n}} \\
&\le& c ~\delta_D(x)^{\alpha/2} \sum^{\infty}_{k=1} \sigma (\partial D \cap B (z_x, 2^k \delta_D (x)) ) 
~(2^{k-2} \delta_D (x) )^{-n}\\
&\le& c ~c_1 \delta_D(x)^{\alpha/2} \left(\sum^{\infty}_{k=1} 2^{2n-k} \right)              \delta_D(x)^{-1}
~\le~ c \delta_D(x)^{\alpha/2-1}  ~~\mbox{ for every } x \in D.
\end{eqnarray*}
\qed

\medskip

If $u$ and $h$ are singular $(-\Delta)^{\alpha/2}$-harmonic functions in $D$ and $u/h$ is bounded, then $u$ can be recovered from non-tangential boundary limit values of $u/h$.

\medskip

\begin{thm}\label{rep}
If $u$ is a singular $(-\Delta)^{\alpha/2}$-harmonic function in $D$ and $u/h$ is bounded for a positive singular $(-\Delta)^{\alpha/2}$-harmonic function $h$ in $D$ with the Martin-representing measure $\nu$, then for every $x \in D$
$$
u(x)~=~h(x)\, \E^h_x\left[\varphi_u \left(\lim_{t \uparrow \tau^h_D} X^h_t \right)\right]
$$
where $$
\varphi_u(z):= \lim_{ A^{\beta}_z \ni x \rightarrow z} \frac{u(x)}{h(x)},~~\beta >\frac{1-\kappa}{\kappa},
$$ which is well-defined for $\nu$-a.e. $z \in \partial D$. If we further assume that $u$ is   positive in $D$, then 
$$
u(x)= \int_{\partial D} M_D(x,w)\, \varphi_u(w)\, \nu(dw)
$$
That is, $\varphi_u(z)$ is 
Radon-Nikodym derivative of the (unique) Martin-representing measure $\mu_u$ with respect to $\nu$.
\end{thm}

\pf
Without loss of generality, we can assume $u$ is positive and bounded.
Recall that
$$\P_{x_0}^{z}\left(\lim_{t\uparrow \tau^z_D} X_t =z\right)\,=\,1$$
for every $z \in \partial D$.
Now take an increasing sequence of smooth open sets $\{D_m\}_{m \ge 1}$ such that $\overline{D_m} \subset D_{m+1}$ and $\cup^{\infty}_{m=1} D_m = D$. Then
\begin{eqnarray*}
1&=&\P_{x_0}^{z}\left(\lim_{m \rightarrow \infty} 
\left(\frac{u}{v}\right)\left(X^z_{\tau^z_{D_m}}\right)
= \lim_{t\uparrow \tau^z_D} \frac{u(X^z_t)}{h(X^z_t)}= \lim_{ A^{\beta}_z \ni x \rightarrow z} \frac{u(x)}{h(x)} \right)\\
&=&\P_{x_0}^{z}\left(\lim_{m \rightarrow \infty} 
\left(\frac{u}{v}\right)\left(X^z_{\tau^z_{D_m}}\right) = \varphi_u(z),~ \lim_{t\uparrow \tau^z_D} X^z_t =z \right)\\
&=&\P_{x_0}^{z}\left(\lim_{m \rightarrow \infty} 
\left(\frac{u}{v}\right)\left(X^z_{\tau^z_{D_m}}\right)
 = \varphi_u\left(\lim_{t\uparrow \tau^z_D} X^z_t\right)\right)
\end{eqnarray*}
for $\nu$-a.e. $z \in \partial D$.
By Proposition \ref{prop:3.3} and zero-one law (Proposition \ref{prop:3.5})
$$
\lim_{m \rightarrow \infty}
\left(\frac{u}{v}\right)\left(X^z_{\tau^z_{D_m}}\right)
\,=\,
 \varphi_u \left( \lim_{t\uparrow \tau^h_D} X^h_t \right) 
~~~\P^h_x \mbox{-a.s. for every} x \in D.
$$
Therefore, by the bounded convergence theorem and the harmonicity of $u/h$ with respect to $\P^h_x$, we have
$$
\frac{u(x)}{h(x)}
\,=\, 
\lim_{m \rightarrow \infty}
 \E^h_x \left[
\left(\frac{u}{v}\right)\left(X^z_{\tau^z_{D_m}}\right)
\right]
 \,=\,\
E^h_x\left[\lim_{m \rightarrow \infty} 
\left(\frac{u}{v}\right)\left(X^z_{\tau^z_{D_m}}\right)
\right]
\,=\,
\E^h_x\left[\varphi_u\left(\lim_{t\uparrow \tau^h_D} X^h_t\right)\right] $$
for every $x \in D.$
By Theorem \ref{thm:son2}, 
$$
u(x)= \int_{\partial D} M_D(x,w) \,\varphi_u(w) \,\nu(dw).
$$

\qed

\medskip

\begin{remark}\label{R1}
{\rm The above theorem is not true without the boundedness assumption. For example, fix a boundary point $z_1 \in \partial D$  and let $u(x):= M_D(x,z_1)$. We let $N$ be the set of irregular boundary points of $D$ and choose any finite measure $\nu$ on $\partial D$ with $\nu(N \cup \{z_1\})=0.$ Then $u$
 can not have that representation in Theorem \ref{rep} with the positive singular $(-\Delta)^{\alpha/2}$-harmonic function $h$ in $D$ with the Martin-representing measure $\nu$.
 In fact, we know from Theorem \ref{T:M} and its proof that for every regular boundary point $z \not= z_1$,
$$
\varphi_u(z)\,=\,\lim_{ A^{\beta}_z \ni x \rightarrow z} \frac{u(x)}{h(x)}
\,=\,\lim_{ A^{\beta}_z \ni x \rightarrow z}u(x) \times
\lim_{ A^{\beta}_z \ni x \rightarrow z} \frac1{h(x)}=0.
$$
Therefore 
$$\int_{\partial D} M_D(x,w) \varphi_u(w) \nu(dw)=0 ~~\mbox{ for every } x \in D,$$
which is obviously not equal to $M_D(x,z_1).$}
\end{remark}

\medskip
Recall that an open set $G \in \RR^n$ is said to satisfy the uniform exterior cone condition if there exist constant 
$\eta >0$, $r >0$ and a cone 
$\C=\{ x =(x_1, \cdots, x_n) \in \RR^n ; 0 < x_n, (x_1^2 + \cdots + x_{n-1}^2)^{1/2} < \eta x_n \}$
such that for every $z \in \partial G$, there is a cone $\C_z$ with vertex $z$, isometric to $\C$ and satisfying 
 $\C_z \cap B(z,r) \subset G^c$. It is well known that for every bounded open set $G$ in $\RR^n$ satisfying uniform exterior cone condition, there is a function $P_G(x,z)$ defined on $G\times G^c$ such that
$$
\E_x[f(X_{\tau_G})] \,=\, \int_{G^c} P_G(x,z)f(z)dz, ~~~ x \in G
$$
for every $f \ge 0$ on $G^c$ (for example, see \cite{CS3}). The kernel   $P_G(x,z)$ is called the Poisson kernel for the symmetric $\alpha$-stable process in $G$.
If our $D$ satisfies the uniform exterior cone condition, then we can state Theorem \ref{rep} for positive $(-\Delta)^{\alpha/2}$-harmonic function in $D$. Note that if  $D$ satisfies the exterior cone condition and $u$ is a positive $(-\Delta)^{\alpha/2}$-harmonic function in $D$, then $u(x)-\E_x[u(X_{\tau_D})]$ is a positive singular $(-\Delta)^{\alpha/2}$-harmonic function (see the proof of Theorem 4.3 in \cite{CS3}).
\medskip

\begin{cor}\label{CC1}
Suppose $D$ satisfies the exterior cone condition.
If $u$ is a positive $(-\Delta)^{\alpha/2}$-harmonic function in $D$ and 
$$\frac{u(x)-\E_x[u(X_{\tau_D})]}{h(x)}
$$ 
is bounded in $D$ for a positive singular $(-\Delta)^{\alpha/2}$-harmonic function $h$ in $D$, then for every $x \in D$
$$
u(x)= \int_{D^c} P_D(x,y) u(y) dy + \int_{\partial D} M_D(x,w) \varphi_u(w) \nu(dw)
$$
where 
$P_D(x,y)$ is Poisson kernel for $X$ in $D$, $\nu$ is the Martin-representing measure for $h$, and
$$
\varphi_u(z):= \lim_{ A^{\beta}_z \ni x \rightarrow z} \frac{u(x)}{h(x)}, \qquad \beta >\frac{1-\kappa}{\kappa},
$$ which is well-defined for $\nu$-a.e. $z \in \partial D$. 
\end{cor}

\medskip
If $D$ is a bounded $C^{1,1}$-open set,
we have a concrete sufficient condition for the representation theorem.
\medskip

\begin{thm}\label{AC}
Suppose $D$ is bounded $C^{1,1}$-open set. Let $\sigma$ be the surface measure on $\partial D$ and let
$$
h(x) := \int_{\partial D} M_D(x, w) ~\sigma (dw), \qquad x \in D.
$$
 If $u$ is positive $(-\Delta)^{\alpha/2}$-harmonic in $D$ and there exists $c >0$ such that
$u(x)-\E_x[u(X_{\tau_D})] \le c \delta_D(x)^{\alpha/2-1} $ for $x \in D$, then 
for every $x \in D$
$$
u(x)= \int_{D^c} P_D(x,y) u(y) dy + \int_{\partial D} M_D(x,w) \varphi_u(w) \sigma(dw),
$$
where 
$P_D(x,y)$ is Poisson kernel for $X$ in $D$ and
$$
\varphi_u(z):= \lim_{ A^{\beta}_z \ni x \rightarrow z} \frac{u(x)}{h(x)}, \qquad \beta > 1,
$$ which is well-defined for $\nu$-a.e. $z \in \partial D$. 
\end{thm}

\pf
It is easy to see that bounded $C^{1,1}$-open sets satisfy the uniform exterior cone condition. 
Let $\sigma$ be the  surface measure on $\partial D$. Since $D$ is a bounded $C^{1,1}$ open set, there exists $c_1$ depending only on $D$ such that for every $x \in D$ and $k \ge 1$,
\begin{equation}\label{A1}
c_1 (2 \delta_D (x))^{n-1} \le~\sigma ( \partial D \cap B (z_x, 2 \delta_D (x)) ),
\end{equation}
where $z_x \in \partial D$ and $|z_x - x| = \delta_D (x)$.
Let 
$$
E := \left\{ w \in \partial D ;~ |w - z_x| \le 2 \delta_D (x) \right\} .
$$
We see that
\begin{equation}\label{A2}
\delta_D(x)~ \le ~|x - w|~ \le~ 3 \delta_D(x) ~~~\mbox{ if } w \in E.
\end{equation}
By the Martin kernel estimate in \cite{C1}, (\ref{A1}) and (\ref{A2}) we have  
$$
h(x) ~\ge~ \delta_D(x)^{\alpha/2} \int_{E} \frac{\sigma(dw)}{|x-w|^{n} }
~\ge~ c  \delta_D (x)^{\alpha/2-n} \sigma(E) 
~\ge~ c c_1 \delta_D (x)^{\alpha/2-1}~~~~ \mbox{ for every } x \in D.
$$
So 
$$\left|\frac{u(x)-\E_x[u(X_{\tau_D})]}{h(x)}\right| ~\le~ C~ <~ \infty~~~~~\mbox{ for every } x \in D.
$$ 
Therefore by Corollary \ref{CC1}, 
$$
u(x)~=~ \int_{D} P_D(x,y) u(y) dy \,+\, \int_{\partial D} M_D(x,w) \varphi_u(w) \sigma(dw).
$$

\qed
\medskip

\medskip
Now suppose that $n=2$, $D=B:=B(0,1)$, $x_0=0$ and $\sigma_1$ is the normalized surface measure on $\partial B$. It is showed in \cite{K} that the Stolz domain is the best possible one for Fatou's theorem in $B$ for transient censored stable processes. Using similar  methods, we can show that our Stolz open set is also the best possible one. First we modify the proof of Lemma 3.19 in \cite{K} using the Martin kernel estimate for symmetric $\alpha$-stable process in $B$. 

\medskip
\begin{lemma}\label{counter_a}
Let 
$$h(x):=\int_{ \partial B } M_B(x, w) \sigma_1(dw).$$
Suppose $U$ is a measurable function on $\partial B$ such that $0 \le U \le 1.$ Let 
$$u(x):=\int_{ \partial B } M_B(x, w) U(w) \sigma_1(dw) = \frac1{2 \pi} \int^{2 \pi}_0 M_B(x, e^{i \theta}) U( e^{i \theta}) d  \theta$$
where $x \in B.$
Suppose that  $0 < \lambda < \pi $ and $U(e^{i \theta})=1$ for $\theta_0 - \lambda \le \theta \le  \theta_0 + \lambda.$
Then there exists a $\delta = \delta (\eps, \alpha)$ such that 
$$
1-\eps \le \frac{u(\rho e^{i \theta_0})}{h(\rho e^{i \theta_0})} \le 1  ~~~~~\mbox{if } \rho > 1- \lambda \delta.
$$ 
\end{lemma}
\pf
First, it is clear that 
$$ \frac{u(x)}{h(x)}=\frac1{h(x)}\int_{ \partial B } M_B(x, w) U(w) \sigma_1(dw) \le \frac1{h(x)}\int_{ \partial B } M_B(x, w) \sigma_1(dw) \equiv 1$$ for every $x \in B.$

Let $V:=\frac1{2}(U-1)$ so that $|V| \le 1$ and $V=0$ for  $\theta_0 - \lambda \le \theta \le  \theta_0 + \lambda.$
If $\frac{1- \rho}{\lambda} < \delta < \frac2{\pi},$
\begin{eqnarray*}
&&|\rho  e^{i \theta_0}- e^{i \theta}| \,\geq\,
  |  e^{i \theta_0}- e^{i \theta}|-(1-\rho)\,
 \geq \,
 2 \left|\sin\left(\frac {\theta_0- \theta }{2}\right)\right|-
 \delta |\theta_0- \theta|  \\ 
&&\geq  \frac2{\pi}|\theta_0- \theta| - \delta |\theta_0- \theta|
\,=\, ( \frac2{\pi}-\delta )|\theta_0- \theta| \qquad \mbox{for } |\theta_0- \theta| > \lambda.
\end{eqnarray*}
So we have, by  the Martin kernel estimate for ball,
\begin{eqnarray*}
\left|\frac1{2 \pi} \int^{2 \pi}_0  M_D(\rho e^{i \theta_0}, e^{i \theta}) V( e^{i \theta}) d  \theta \right| 
&\le & c \int^{2 \pi}_0 \frac{|V( e^{i \theta})|}{|\rho  e^{i \theta_0}- e^{i \theta}|^{2}  }
  d  \theta~ (1-\rho)^{\alpha/2}\\
 &\le & c (1-\rho)^{\alpha/2} ( \frac2{\pi}-\delta )^{-2} \int_{|\theta_0- \theta| > \lambda} \frac{d \theta }{|\theta_0- \theta|^{2} }\\
 &\le & c (1-\rho)^{\alpha/2} ( \frac2{\pi}-\delta )^{-2} \lambda^{-1}\\
 &\le &c \delta(\frac2{\pi} - \delta)^{-2} (1-\rho)^{\alpha/2-1}.
\end{eqnarray*}
Therefore, if $ \delta \le \frac1{\pi},$
\begin{eqnarray*}
\frac{ u(\rho e^{i \theta_0}) }{h(\rho e^{i \theta_0})}&=&\frac1{h(\rho e^{i \theta_0})} \frac1{2 \pi} \int^{2 \pi}_0  M_D(\rho e^{i \theta_0}, e^{i \theta}) (1+2V( e^{i \theta})) d  \theta\\
&\geq&\frac1{h(\rho e^{i \theta_0})} \left(h(\rho e^{i \theta_0})-c \delta(\frac2{\pi} - \delta)^{-2} (1-\rho)^{\alpha/2-1}\right)
~\geq~1-c_1 \delta.
\end{eqnarray*}
In the above inequality, we have used the fact that
$h(\rho e^{i \theta_0}) \le c (1-\rho)^{\alpha/2-1}$ (see the proof of Corollary \ref{C2:Fatou}).
For any $\eps >0$, $\delta:=\min\left\{\frac{\eps}{c_1}, \frac1{\pi} \right\}$ will do.
\qed

\medskip

Once we have this lemma, the rest of the details are similar to those in \cite{L}.
A curve $C_0$ is called a tangential curve in $B$ which ends on $\partial B$ if $C_0 \cap \partial B = \{w_0\} \in \partial B$, $C_0 \setminus \{w_0\} \subset B$ and there are no $r > 0$ and $ \beta >1$ such that $C_0 \cap B(w_0,r) \subset A^{\beta}_{w_0} \cap B(w_0,r)$.  

\medskip

\begin{thm}\label{counter_thm}
Let 
$$h(x):=\int_{ \partial B } M_B(x, w) \sigma_1(dw).$$ Let $C_0$ be a tangential curve in $B$ which ends on $\partial B$ and let $C_{\theta}$ be the rotation of $C_{0}$ about $x_0$ through an angle $\theta.$
Then there exists a positive $(-\Delta)^{\alpha/2}$-harmonic function $u$ in $B:=B(x_0,1)$ such that for a.e. $\theta \in [0, 2\pi]$ with respect to Lebesgue measure,
$$
\lim_{|x| \rightarrow 1, x \in C_{\theta}} \frac{u(x)}{h(x)} \mbox{ does not exist}.
$$
\end{thm}
\pf
We observe that for bounded measurable function $U \ge 0$ on $\partial B$,
$$
\frac1{h(x)}\int_{ \partial B } M_B(x, w) U(w) \sigma_1(dw)
$$
tends radially to $U$ for almost all $\theta$ by the uniqueness of the Martin-representing measure and Theorem \ref{rep}.
So, with the same $E_k$ defined in \cite{L}, one can show that 
there exists a singular $(-\Delta)^{\alpha/2}$-harmonic function $u_k$ satisfying $ 0 \le u_k \le 2^{-k}$ such that, for a set $E^{*}_k$
 equivalent to $E_k$ (i.e. $\sigma_1((E^{*}_k \setminus E_k) \cup (E_k \setminus E^*_k))=0$), 
$$
\lim \frac{u_k}{h}=0 ~\mbox{ radially ~and } ~\limsup    \frac{u_k}{h}= 2^{-k} \mbox{~ along one branch of } C_{\theta}.
$$  
By following the argument in  \cite{L}, one can check that $u:= \sum_{k=1}^{\infty} u_k$ will do.
\qed

\medskip

\section{Relative Fatou's Theorem under Nonlocal Feynman-Kac Transforms}

First we will show the existence of nontangential limit of the ratio of Green function 
and singular $(-\Delta)^{\alpha/2}$-harmonic function. This result will be used later in this section. 
\begin{lemma}\label{New1}
Let $\nu$ be a finite measure on $\partial D$ and let
$$
h(x) := \int_{\partial D} M_D(x, w) \,\nu (dw),~~~ x \in D.
$$ 
then for $\nu$-a.e. $z \in \partial D$ and every $y \in D$,
$$
\lim_{ A^{\beta}_z \ni x \rightarrow z} \frac{G_D(x, y)}{h(x)} \mbox{ exists for every } \beta >\frac{1-\kappa}{\kappa}.
$$
The above nontangential limit typically depends on $y \in D$ but the null set for the limit is independent of $y \in D$.
\end{lemma}
\pf
It is well known that (for example, see \cite{B}, \cite{CS3} and \cite{CZ})
$$
G_D(x,y)=G(x,y)-\E_x[G(X_{\tau_D},y)]
$$
 where
$G(x,y)$ is the Green function of $X$ in $\RR^n$.
Since $G(x,x_0)$ is continuous near $\partial D$,  for every $z \in \partial D$
$$
\lim_{ A^{\beta}_z \ni x \rightarrow z} G(x, x_0)~=~G(z,x_0)~~~ \mbox{ for every } \beta >\frac{1-\kappa}{\kappa}.
$$
Therefore by Theorem \ref{T:Fatou}, for $\nu$-a.e. $z \in \partial D$,
$$
\lim_{ A^{\beta}_z \ni x \rightarrow z} \frac{G_D(x, x_0)}{h(x)} 
=\lim_{ A^{\beta}_z \ni x \rightarrow z} \frac{1}{h(x)}\times \lim_{ A^{\beta}_z \ni x \rightarrow z} G(x, x_0)-
\lim_{ A^{\beta}_z \ni x \rightarrow z}  \frac{\E_x[G(X_{\tau_D},x_0)]}{h(x)} 
$$
exists for every $\beta > (1-\kappa)/\kappa.
$
Since
 $$
\lim_{x\rightarrow z \in \partial D} \frac{G_D(x, y)}{G_D(x, x_0)} =\lim_{x\rightarrow z \in \partial D} \frac{G_D(y, x)}{G_D(x_0,x)}=M_D(y,z)
$$
 for every
$z\in \partial D$ and  $y \in D$, we have  for $\nu$-a.e. $z \in \partial D$ and every $y \in D$,
$$
\lim_{ A^{\beta}_z \ni x \rightarrow z} \frac{G_D(x, y)}{h(x)} = \lim_{ A^{\beta}_z \ni x \rightarrow z} \frac{G_D(x, y)}{G_D(x,x_0)} \times
\lim_{ A^{\beta}_z \ni x \rightarrow z} \frac{G_D(x, x_0)}{h(x)}~ \mbox{ exists for every } \beta >\frac{1-\kappa}{\kappa}.
$$
\qed

\medskip

Now we recall the following definitions from Chen
\cite{C} and specify them
for $X^D$, the symmetric $\alpha$-stable process in $D$.
We call a positive measure $\mu$ on $D$ a smooth measure
of $X^D$ if there is a positive continuous additive functional
(PCAF in abbreviation) $A$ of $X^D$ such that
\begin{equation}\label{eqn:Revuz1}
\int_D f(x) \mu (dx) = \uparrow \lim_{t\downarrow 0}
\int_D \E_x \left[ \frac1t \int_0^t f(X^D_s) dA_s \right] dx
\end{equation}
for any Borel measurable function $f\geq 0$.
Here $\uparrow \lim_{t\downarrow 0}$ means the quantity
is increasing as $t\downarrow 0$.
The measure $\mu$ is called the Revuz measure of $A$.
It is known that $\E_x [A_{\tau_D} ] = \int_D G_D(x, y) \mu (dy)$.
For a signed measure $\mu$, we use $\mu^+$ and $\mu^-$
to denote its positive and negative parts respectively.
If $\,\mu^+$and $\mu^{-}$ are smooth measures of $X^D$ and
$A^{+}$ and $A^{-}$ are their corresponding PCAFs
of $X^D$, then we say the continuous additive functional 
$A:=A^+-A^{-}$ of $X^D$ has (signed) Revuz measure $\mu$.
Let $d$ denote the diagonal of $D\times D$.

\medskip

\begin{defn}\label{df:5.1}
Suppose that $A$ is a continuous additive
functional of $X^D$ with Revuz measure $\nu$. Let
$A^{+}$ and $A^{-}$ be the PCAFs (positive continuous additive
functionals) of $X^D$ with Revuz measures $\nu^{+}$ and $\nu^{-}$ respectively.
Let $|A|= A^{+}+A^{-}$ and $|\nu|=\nu^{+}+\nu^{-}$.

\begin{description}
\item{(1)}
The measure $\nu$ (or the continuous additive functional $A$)
is said to be in the class $\S_\infty (X^D)$
if for any $\eps>0$ there is a Borel
subset $K=K(\eps)$ of finite $| \nu|$-measure
and a constant $\delta = \delta (\eps) >0 $ such that
$$ 
\sup_{(x, z)\in (D\times D) \setminus d} \int_{D\setminus K}
\frac{G_D(x, y)G_D(y, z) }{ G_D(x, z)} \, |\nu | (dy)\le\eps
$$
and for all measurable set $B\subset K$ with $| \nu |(B)<\delta$,
$$
 \sup_{(x, z)\in (D\times D) \setminus d}
 \int_B \frac{G_D(x, y)G_D(y, z) }{ G_D(x, z)} \, |\nu |
  (dy) \le\eps.
$$

\item{(2)}
A function $q$ is said to be in the class
$\S_\infty (X^D)$,
 if $\nu(dx) :=q(x)dx$ is in the class
$\S_\infty (X^D)$.
\end{description}
\end{defn}
\medskip
Let $(N^D, H^D)$ be a L\'evy system for the symmetric  
$\alpha$-stable process $X^D$ in $D$, 
that is, for every $x\in D$, $N^D(x, dy)$  is a kernel on
$(D_\partial, {\cal B}(D_\partial))$, where $\partial $ is 
the  cemetery point for process $X^D$ 
and $D_\partial = D\cup \{ \partial \}$, 
and $H^D_t$ is a positive continuous additive 
functional of $Y$ with bounded
1-potential such that  for any nonnegative Borel function $f$ on 
$D \times D_\partial $
that vanishes along the diagonal $d$, 
$$
\E_x\left(\sum_{s\le t}f(X^D_{s-}, X^D_s) \right)
= \E_x\left(\int^t_0\int_{D_\partial}
 f(X^D_s, y) N(X^D_s, dy)dH^D_s\right)
$$
for every $x\in D$ (see \cite{Sharpe} for details).
We let $\mu_{H^D}(dx)$ be the Revuz measure for $H^D$.
\medskip

\begin{defn}\label{df:5.2} 
Suppose $F$ is a bounded function
on $D\times D$ vanishing on the diagonal. Let
$$
\mu_{|F|} (dx):=\left( \int_D |F(x, y)| N^D(x, dy)\right) \mu_{H^D}(dx).
$$
 $F$ is said to be in the class $\A_\infty(X^D)$
if for any $\eps>0$ there is a Borel
subset $K=K(\eps)$ of finite $\mu_{|F|}$-measure
and a constant $\delta = \delta (\eps) >0 $ such that
$$
\sup_{(x, w)\in (D\times D) \setminus d} \int_{(D\times D) \setminus (K\times K)} G_D(x, y)
\frac{|F(y, z)|G_D(z, w) }{ G_D(x, w)}N^D(y,dz)\mu_{H^D}(dy) \,\le\,\eps
$$
and for all measurable sets $B\subset K$ with $\mu_{|F|} (B)<\delta$,
$$
\sup_{(x, w)\in (D\times D) \setminus d}
\int_{(B\times D)\cup (D\times B)}
G_D(x, y)
\frac{|F(y, z)| G_D(z, w) }{G_D(x, w)}N^D(y,dz)\mu_{H^D}(dy)\,\le\,\eps.
$$
\end{defn}

\medskip

As it is remarked in Chen \cite{C}, it follows from measure
theory that the Borel set in above Definitions 
\ref{df:5.1}-\ref{df:5.2} can be taken to be compact.

For a smooth measure $\mu$ associated with a continuous additive
functional $A^\mu$ and a Borel measurable function $F$ on $D\times D$ that
vanishes along the diagonal, define
$$
e_{A^\mu+ F}(t):=\exp \left( A^\mu_t + \sum_{0<s\le t}F(X^D_{s-}, X^D_s)
\right), \quad t\geq 0.
$$
In the remainder of this section,
let $\mu\in \S_\infty(X^D)$ and $F\in \A_\infty (X^D)$
such that the gauge function $x\mapsto \E_x \left[
e_{A^\mu+F}(\tau_D)\right]$ is bounded.
It leads us a Schr\"odinger semigroup
$$
 Q_t f(x):= \E_x \left[ e_{A^\mu+F}(t)f(X^D_t) \right], ~~~~x \in D.
$$
For $x, y \in D$, let $\E_x^y$ denote the expectation 
for the conditional process starting from $x$ 
obtained from $X^D$ through
Doob's $h$-transform with $h(\,\cdot\,)= G_D(\,\cdot\, , y)$.
By Lemma 3.9 of Chen \cite{C}, the Green function for the
Schr\"odinger semigroup $\{Q_t, \, t\geq 0\}$  is
$V_D(x, y)=
   u(x,y)   G_D(x, y)$ where $u(x,y):= \E_x^y\left[ e_{A^\mu+F}(\tau_D^y) \right]$,
that is,
$$
\int_D V_D(x, y) f(y) \, dy = \int_0^\infty Q_t f (x) \, dt
= \E_x \left[ \int_0^\infty e_{A^\mu +F}(t) f(X^D_t) \, dt \right]
$$
for any Borel measurable function $f\geq 0$ on $D$. Thus
$V_D(x, y)$ is comparable to $G_D(x, y)$ on $(D\times D) \setminus d$ by Theorem 3.10 in \cite{C}.

For $x \in D$ and $ w \in \partial D$, let $u(x, w):=\E_x^w \left[ e_{A^\mu+F}(\tau_D^w) \right]$ where $\E_x^w$ is the expectation for
the conditional process of $X^D$ obtained through $h$-transform
with $h(\,\cdot \,)=M_D(\,\cdot\, , w)$. One can follow the argument in Section 3 of \cite{CK2} and show that, for any $w \in  \partial D$ and  $x \in D$,
\begin{equation}\label{UMM}
u(x, w)= \lim_{  D \ni y \to w} \E_x^y \left[ e_{A^\mu+F}(\tau_D^y) \right]
\end{equation}
and 
\begin{eqnarray}
(u(x, w)-1)M_D(x, w)=\int_D V_D(x, z)  M_D(z, w) \mu(dz) ~~~~~~~~~~~~~~~~~~~~~~~~~~~~~~~~~~~~~~~~~~~~~~\nonumber\\
~~~~~~~~~~~~~~~~~~+\int_D V_D (x, y)
\left( \int_D (e^{F(y, z)}-1)
M_D(z, w)N^D(y, dz)\right) \mu_{H^D}(dy), \label{UM} 
\end{eqnarray}
which implies that for every $w \in  \partial D$ and  $x \in D$,
\begin{equation}\label{KM}
K_D(x, w):= \lim_{ D \ni y \to w} \frac{V_D(x, y)}{V_D(x_0, y)}=M_D(x, w) \,  \frac{u(x, w)}{u(x_0, w)} \approx M_D(x, w).
\end{equation}
 The above is proved for transient censored stable process in bounded $C^{1,1}$-open set in Section 3 of \cite{CK2}. However,  the same proof works for $X^D$ in bounded $\kappa$-fat open set $D$.

\medskip

\begin{lemma}\label{New2}
Let $\nu$ be a finite measure on $\partial D$ and let
$$
h(x) := \int_{\partial D} M_D(x, w) \,\nu (dw),~~~ x \in D.
$$ 
Then for $\nu$-a.e. $z \in \partial D$ and every $y \in D$,
\begin{equation}\label{AS}
\lim_{ A^{\beta}_z \ni x \rightarrow z} \frac{V_D(x, y)}{h(x)} \mbox{ exists for every } \beta >\frac{1-\kappa}{\kappa}.
\end{equation}
Moreover, for $\nu$-a.e. $z \in \partial D$
\begin{equation}\label{AAS}
\lim_{ A^{\beta}_z \ni x \rightarrow z}  \frac1{h(x)}\left[\int_D V_D(x, y) f(y)\mu(dz)+ \int_{D \times D} V_D (x, y)
( e^{F(y, z)}-1 )f(z)
 N^D(y, dz) \mu_{H^D}(dy) \right]~ 
\end{equation}
 exists for every $(-\Delta)^{\alpha/2}$-harmonic function $f \ge 0$ and every $ \beta >(1-\kappa)/\kappa$.
\end{lemma}

\pf
 Note that 
the definition of $\A_\infty (X^D)$ is symmetric in $x$ and $y$
and so $\wh F\in \A_\infty(X^D)$ where $\wh F(x, y):=F(y, x)$.
Moreover by the argument in page 60 of \cite{CS7} and 
the symmetric property of $X^D$, it is easy to see that
$$
\E_x^y \left[e_{A^\mu+F}(\tau_D^y) \right] 
\,= \,\E_y^x \left[e_{A^\mu+ \wh F}(\tau_D^x) \right]
~~~\mbox{ for } x, y\in D,
$$
which also implies that $\E_y [e_{A^\mu+ \wh F}(\tau_D) ]$ is bounded by Theorem 3.10 in \cite{C}.
Thus by (\ref{UMM}) with $\wh F$ instead of $F$, we have for every $ y\in D$ and $w \in \partial D$,
$$
\lim_{D \ni x \to w} u(x,y)\,=\,\lim_{ D \ni x \to w}  \E_y^x \left[e_{A^\mu+ \wh F}(\tau_D^x) \right]
\,=\, \E_y^w \left[e_{A^\mu+ \wh F}(\tau_D^w) \right].
$$
So (\ref{AS}) is true by Lemma \ref{New1}.

Now we will show (\ref{AAS}) with a fixed $(-\Delta)^{\alpha/2}$-harmonic function $f \ge 0$.
By Theorem \ref{T:Fatou} and (\ref{AS}), 
there is a $\nu$-null set $ B \subset \partial D$ such that for every $ w \in \partial D \setminus B $ , 
$$
\lim_{A^{\beta}_w \ni x \rightarrow w} \frac{f(x)}{h(x)},~\lim_{A^{\beta}_w \ni x \rightarrow w} \frac1{h(x)}~\mbox{ and }~\lim_{A^{\beta}_w \ni x \rightarrow w} \frac{V_D(x, y)}{h(x)}~ \mbox{ exist for every } \beta >\frac{1-\kappa}{\kappa},  y\in D.
$$
Now we fix $ \beta >(1-\kappa)/\kappa$, $ w_0 \in \partial D \setminus B$ and 
a sequence $\{x_k\}_{k \ge 1} \subset  A^{\beta}_{w_0} $ converging to $ w_0$.
Let
$$ M:= \sup_{(x, y) \in D \times D \setminus d} \E_x^y \left[ e_{A^\mu+F}(\tau_D^y) \right] < \infty
~\mbox{ and }~ \wt F (y,z):= e^{F(y, z)}-1
$$
so that $e^{|F(y, z)|}-1 \ge |\wt F (y,z)|=\wt F^+ (y,z)+\wt F^- (y,z)$ where $\wt F^+ (y,z)$ and $\wt F^- (y,z)$ are $\wt F$'s positive and negative parts respectively.
Given $\eps > 0 ,$ by a similar argument to those for Proposition 3.1 in \cite{CS7} and Proposition 3.1 in \cite{CS6} (also see the remark immediately following Definition \ref{df:5.2}), there exists a compact subset $K = K(\eps, M) \subset D$ such that
$$
\int_{D \setminus K} G_D(x,y) f(y) |\mu|(dy)
~ \le~ \frac{\eps f(x)} {2M}~~~ \mbox{ for all } x \in D 
$$ 
and
$$
\int_{(D \times D) \setminus (K \times K)} G_D(x,y) ( e^{|F(y,z)|} - 1) f(z) N^D (y,dz) \mu_{H^D}(dy) 
~\le~ \frac{\eps f(x)}{2M}
$$
for all $x \in D.$
Thus, for every $x \in D$, 
\begin{eqnarray*}
 &&\frac1{h(x)}\left[\int_D V_D(x, y) f(y)\mu^+(dz)+ \int_{D \times D} V_D (x, y)
   \wt F^+ (y,z) f(z)
 N^D(y, dz) \mu_{H^D}(dy) \right]\\
&\le&  \frac1{h(x)} \left[ \int_K V_D(x, z) f(z) \mu^+(dz)+ M\int_{D \setminus K} G_D(x, z) f(z) |\mu|(dz) \right.\\
&&~~~~~~+ \int_{K \times K} V_D (x, y)  \wt F^+ (y,z) f(z) N^D(y, dz)  \mu_{H^D}(dy) \\
&&~~~~~~+\left. M \int_{(D \times D) \setminus (K \times K) } G_D (x, y)
( e^{|F(y, z)|}-1 )
f(z) N^D(y, dz)  \mu_{H^D}(dy)  \right]\\
& \le & \frac{\eps f(x)}{h(x)} + \int_K \frac{V_D(x, z)}{h(x)}f(z) \mu^+(dz)+ \int_{K \times K} \frac{V_D (x, y)}{h(x)}
F^+ (y,z)f(z)
 N^D(y, dz) \mu_{H^D}(dy) .
\end{eqnarray*}
Since $1_D$ is a excessive function for $X^D$, by Proposition 3.1 in \cite{CS7} and Proposition 3.1 in \cite{CS6}, we can easily see that
$$
 \{ G_D(x,z) |\mu|(dz) \,;\, x \in D \}
$$
 is uniformly integrable in $D$ and 
$$
 \left\{ G_D(x,y)  ( e^{|F(y, z)|}-1 )
 N^D(y, dz) \mu_{H^D}(dy);~ x \in D \right\}
$$
  is uniformly integrable in $D \times D$.
Thus, since $V_D(x, y)$ is comparable to $G_D(x, y)$ on $(D\times D) \setminus d$, 
$$
 \left\{ \frac{V_D(x_k,z)} {h(x_k)} f(z) 1_K(z) \mu^+(dz) \,; \,k \ge 1 \right\}
$$
 is uniformly integrable in $D$ and 
$$
 \left\{ \frac{V_D(x_k,y)} {h(x_k)} \wt F^+ (y,z) f(z)1_{K \times K} (y,z)
 N^D(y, dz) \mu_{H^D}(dy)\,;\, k \ge 1 \right\}
$$
  is uniformly integrable in $D \times D$ because $f$ is bounded on $K \times K$ and 
$
\lim_{k \to \infty} \frac1{h(x_k)}
$ exists.
Therefore, 
\begin{eqnarray*}
&&\limsup_{k \to \infty} \frac1{h(x_k)}\left[\int_D V_D(x_k, z) f(z)\mu^+(dz)+ \int_{D \times D} V_D (x_k, y)
   \wt F^+ (y,z) f(z)
 N^D(y, dz) \mu_{H^D}(dy) \right]\\ 
&\le&  \eps \lim_{k \to \infty}\frac{f(x_k)}{h(x_k)} ~+~ \int_K \lim_{k \to \infty}\frac{V_D(x_k, z)}{h(x_k)}f(z) \mu^+(dz)\\
&&~~~~~~~~~~~~~~~~~~~~~~~~~~~~~~+ \int_{K \times K} \lim_{k \to \infty}\frac{V_D (x_k, y)}{h(x_k)}
\wt F^+ (y,z)f(z)
 N^D(y, dz) \mu_{H^D}(dy) \\
&\le&  \eps \lim_{k \to \infty}\frac{f(x_k)}{h(x_k)}~+~ \int_D \lim_{k \to \infty}\frac{V_D(x_k, z)}{h(x_k)}f(z) \mu^+(dz)\\
&&~~~~~~~~~~~~~~~~~~~~~~~~~~~~~~+ \int_{D \times D} \lim_{k \to \infty}\frac{V_D (x_k, y)}{h(x_k)}
\wt F^+ (y,z)f(z)
 N^D(y, dz) \mu_{H^D}(dy) .
\end{eqnarray*}
Letting $\eps \to 0$, we have
\begin{eqnarray*}
&& \limsup_{k \to \infty} \frac1{h(x_k)}\left[\int_D V_D(x_k, z) f(z)\mu^+(dz)+ \int_{D \times D} V_D (x_k, y)
   \wt F^+ (y,z) f(z)
 N^D(y, dz) \mu_{H^D}(dy) \right]\\ 
& \le &  \int_D \lim_{k \to \infty}\frac{V_D(x_k, z)}{h(x_k)}f(z) \mu^+(dz)+ \int_{D \times D} \lim_{k \to \infty}\frac{V_D (x_k, y)}{h(x_k)}
\wt F^+ (y,z)f(z)
 N^D(y, dz) \mu_{H^D}(dy).
\end{eqnarray*}
On the other hand, by Fatou's lemma,
\begin{eqnarray*}
&& \liminf_{k \to \infty} \frac1{h(x_k)}\left[\int_D V_D(x_k, z) f(z)\mu^+(dz)+ \int_{D \times D} V_D (x_k, y)
   \wt F^+ (y,z) f(z)
 N^D(y, dz) \mu_{H^D}(dy) \right]\\ 
& \ge &  \int_D \lim_{k \to \infty}\frac{V_D(x_k, z)}{h(x_k)}f(z) \mu^+(dz)+ \int_{D \times D} \lim_{k \to \infty}\frac{V_D (x_k, y)}{h(x_k)}
\wt F^+ (y,z)f(z)
 N^D(y, dz) \mu_{H^D}(dy).
\end{eqnarray*}
Therefore,
\begin{eqnarray*}
&& \lim_{k \to \infty} \frac1{h(x_k)}\left[\int_D V_D(x_k, z) f(z)\mu^+(dz)+ \int_{D \times D} V_D (x_k, y)
   \wt F^+ (y,z) f(z)
 N^D(y, dz) \mu_{H^D}(dy) \right]\\ 
& = &  \int_D \lim_{k \to \infty}\frac{V_D(x_k, z)}{h(x_k)}f(z) \mu^+(dz)+ \int_{D \times D} \lim_{k \to \infty}\frac{V_D (x_k, y)}{h(x_k)}
\wt F^+ (y,z)f(z)
 N^D(y, dz) \mu_{H^D}(dy).
\end{eqnarray*}
Similarly, we have
\begin{eqnarray*}
&& \lim_{k \to \infty} \frac1{h(x_k)}\left[\int_D V_D(x_k, z) f(z)\mu^-(dz)+ \int_{D \times D} V_D (x_k, y)
   \wt F^- (y,z) f(z)
 N^D(y, dz) \mu_{H^D}(dy) \right]\\ 
& = &  \int_D \lim_{k \to \infty}\frac{V_D(x_k, z)}{h(x_k)}f(z) \mu^-(dz)+ \int_{D \times D} \lim_{k \to \infty}\frac{V_D (x_k, y)}{h(x_k)}
\wt F^- (y,z)f(z)
 N^D(y, dz) \mu_{H^D}(dy).
\end{eqnarray*}
Consequently  
\begin{eqnarray*}
&& \lim_{k \to \infty} \frac1{h(x_k)}\left[\int_D V_D(x_k, z) f(z)\mu(dz)+ \int_{D \times D} V_D (x_k, y)
   ( e^{F(y, z)}-1 )f(z)
 N^D(y, dz) \mu_{H^D}(dy) \right]\\ 
& = &  \int_D \,\lim_{ A^{\beta}_{w_0} \ni x \rightarrow w_0} \frac{V_D(x, z)}{h(x)}\,f(z) \mu(dz)\\
&&~~~~~~~~~~~~~~~~~~~~~~~~+ \int_{D \times D} \,\lim_{A^{\beta}_{w_0} \ni x \rightarrow w_0} \frac{V_D (x, y)}{h(x)}\,
( e^{F(y, z)}-1 )f(z)
 N^D(y, dz) \mu_{H^D}(dy)
\end{eqnarray*}
for every sequence $\{x_k\}_{k \ge 1} \subset  A^{\beta}_{w_0} $ converging to $ w_0$.
\qed

\medskip

To state relative Fatou's theorem for the Schr\"odinger operator corresponding to $Q_t$, we need the following definition. 

\medskip

\begin{defn}\label{def:har2}A  Borel measurable function
$u$ defined on $D$ is said to be $(\mu, F)$-harmonic
 if 
$$ \E_x \left[ e_{A^\mu+F}(\tau_B) |u(X^D_{\tau_B})| \right]
  <\infty \quad \mbox{and} \quad
 \E_x \left[ e_{A^\mu+F}(\tau_B) u(X^D_{\tau_B})
  \right] = u(x), \qquad  x\in B,
$$
for every open set $B$ whose closure is a compact subset of $D$.
\end{defn}
\medskip

In Chen and Kim \cite{CK2}, an integral representation of nonnegative excessive
functions for the Schr\"odinger operator  is established. Moreover it is shown that the Martin boundary is stable 
under non-local Feynman-Kac perturbation. These results hold for a large class of strong Markov processes.
We state a simpler version with respect to $X^D$ here for later use.
\medskip
\begin{thm}\label{TTTT} (Theorem 5.16 in \cite{CK2}. Also see Section 6 in \cite{CK2} for a general setting) 
For every positive $(\mu , F)$-harmonic function $u$, 
there is a unique finite measure $\nu$ on $\partial D$ such that 
\begin{equation}\label{eqn:5.8}
 u(x) =\int_{\partial D} K_D(x, z) \nu (dz).
\end{equation}
\end{thm}
\medskip

We are now in the position to show relative Fatou's theorem for  $(\mu, F)$-harmonic function. The proof is similar to the proof of Theorem 4.4 in \cite{K} but it requires more works (Lemma \ref{New1} and Lemma \ref{New2}). One can see that the complication comes from the irregularity of the boundary of $D$.   
\medskip

\begin{thm}\label{T:Fatou_FK}
Let $D$ be a bounded $\kappa$-fat open set and $\nu$ be a finite measure on $\partial D$.
Let $k$ be a positive $(\mu , F)$-harmonic function with the Martin-representing measure $\nu$.
That is, 
$$
k(x) = \int_{\partial D} K_D(x, w) \,\nu (dw),~~~ x \in D
$$ 
where $\nu$ is a finite measure on $\partial D$.
If $u$ is a nonnegative $(\mu, F)$-harmonic function, then for $\nu$-a.e. $z \in \partial D$,
$$
\lim_{ A^{\beta}_z \ni x \rightarrow z} \frac{u(x)}{k(x)} \mbox{ exists for every } \beta >\frac{1-\kappa}{\kappa}.
$$
\end{thm}

\pf
For $x \in D$ and $ w \in \partial D$, recall $u(x, w)=\E_x^w \left[ e_{A^\mu+F}(\tau_D^w) \right]$ where $\E_x^w$ is the expectation for
the conditional process of $X^D$ obtained through $h$-transform
with $h(\,\cdot \,)=M_D(\,\cdot\,, w)$. 
By Theorem \ref{TTTT}, there is a finite
measure $\mu_1$ on $\partial D$ such that
$$
u(x) =\int_{\partial D} K_D(x,w) \mu_1 (dw) , ~~~x \in D.
$$
Let 
$$
\overline{\mu_1} (dw):=\frac{\mu_1 (dw)}{u(x_0, w)},
$$ 
which is a finite
measure on $\partial D$ because of (3.16) in \cite{C}.
Using (\ref{UM}) and (\ref{KM}), we have
\begin{eqnarray*}
u(x) &=&\int_{\partial D} K_D(x,w) u(x_0, w) \overline{\mu_1} (dw)\\
&=& \int_{\partial D}  M_D(x, w) u(x, w) \overline{\mu_1} (dw)\\
&=& \int_{\partial D}  M_D(x, w) \overline{\mu_1} (dw)+\int_{\partial D}  M_D(x, w) (u(x, w)-1) \overline{\mu_1} (dw)\\
&=& \int_{\partial D}  M_D(x, w) \overline{\mu_1} (dw)+\int_{\partial D}  \left[ \int_D V_D(x, z)  M_D(z, w) \mu(dz) \right.\\
&&~~~~~~~~~~~~ \left.+ \int_D V_D (x, y)
\left( \int_D (e^{F(y, z)}-1)
M_D(z, w)N^D(y, dz)\right)  \mu_{H^D}(dy)
 \right]\overline{\mu_1} (dw).
\end{eqnarray*}
Let 
$$
f(x):= \int_{\partial D}  M_D(x, w) \overline{\mu_1} (dw)~~\mbox{ and }~~
h(x) := \int_{\partial D} M_D(x, w) ~\nu (dw), ~~~x \in D,
$$
 which are $(-\Delta)^{\alpha/2}$-harmonic in $D$ (and are continuous in $D$) and let 
\begin{eqnarray*}
&& g(x) := \int_{\partial D}  \left[ \int_D V_D(x, z)  M_D(z, w) \mu(dz) \right.\\
&&~~~~~~~~~~~~~~~~~\left.+ \int_D V_D (x, y)
\left( \int_D (e^{F(y, z)}-1)
M_D(z, w)N^D(y, dz)\right)  \mu_{H^D}(dy)
 \right]\overline{\mu_1} (dw).
\end{eqnarray*}
By Tonelli's Theorem, for every $ x \in D$, we have
\begin{eqnarray*}
 &&\int_{\partial D}  \left[ \int_D V_D(x, z)  M_D(z, w) |\mu|(dz)\right.\\
&&~~~~~~~\left.+ \int_D V_D (x, y)
\left( \int_D (e^{|F(y, z)|}-1)
M_D(z, w)N^D(y, dz)\right)  \mu_{H^D}(dy)
 \right]\overline{\mu_1} (dw)\\
&=&\int_D V_D(x, z) f(z)  |\mu|(dz)+ \int_{D \times D} V_D(x, y)
\left( e^{|F(y, z)|}-1 \right)
f(z) N^D(y, dz)  \mu_{H^D}(dy)\\
&\le& M \left[\int_D G_D(x, z) f(z)  |\mu|(dz)+ \int_{D \times D} G_D (x, y)
\left( e^{|F(y, z)|}-1 \right)
f(z) N^D(y, dz)  \mu_{H^D}(dy) \right],
\end{eqnarray*}
where $$ M:= \sup_{(x, y) \in D \times D \setminus d} \E_x^y \left[ e_{A^\mu+F}(\tau_D^y) \right] < \infty.$$

Given $\eps > 0 ,$ using an argument similar to that in Lemma \ref{New2}, there exists a compact subset $K = K(\eps, M) \subset D$ such that
$$
\int_{D \setminus K} G_D(x,y) f(y) |\mu|(dy)
~ \le~ \frac{\eps f(x)} {2M}~~~ \mbox{ for all } x \in D 
$$ 
and
$$
\int_{(D \times D) \setminus (K \times K)} G_D(x,y) ( e^{|F(y,z)|} - 1) f(z) N^D (y,dz) \mu_{H^D}(dy) 
~\le~ \frac{\eps f(x)}{2M}
$$
for all $x \in D.$
Thus, for every $x \in D$, 
\begin{eqnarray*}
&&\int_{\partial D}  \left[ \int_D V_D(x, z)  M_D(z, w) |\mu|(dz)\right.\\
&&~~~~~~\left.+ \int_D V_D (x, y)
\left( \int_D (e^{|F(y, z)|}-1)
M_D(z, w)N^D(y, dz)\right)  \mu_{H^D}(dy)
 \right]\overline{\mu_1} (dw)\\
  &\le&  M \left[ \int_K G_D(x, z) f(z) |\mu|(dz)+ \int_{D \setminus K} G_D(x, z) f(z) |\mu|(dz) \right.\\
&&~~~~~~+ \int_{K \times K} G_D (x, y)( e^{|F(y, z)|}-1 )f(z) N^D(y, dz)  \mu_{H^D}(dy) \\
&&~~~~~~+\left. \int_{(D \times D) \setminus (K \times K) } G_D (x, y)
( e^{|F(y, z)|}-1 )
f(z) N^D(y, dz)  \mu_{H^D}(dy)  \right]\\
& \le & MN \left[\int_K G_D(x, z) |\mu|(dz)+ \int_{K \times K} G_D (x, y)
( e^{|F(y, z)|}-1 )
 N^D(y, dz) \mu_{H^D}(dy) \right]+~\eps f(x) <\infty
\end{eqnarray*}
where 
$$ 
N:= \sup_{y \in K} f(y) < \infty.
$$
So by Fubini's theorem,
$$
g(x)=\left[\int_D V_D(x, y) f(y)\mu(dz)+ \int_{D \times D} V_D (x, y)
( e^{F(y, z)}-1 )f(z)
 N^D(y, dz) \mu_{H^D}(dy) \right].
$$
Therefore by Lemma \ref{New2}, for $\nu$-a.e. $z \in \partial D$
$$
\lim_{ A^{\beta}_z \ni x \rightarrow z} \frac{g(x)}{h(x)}  ~\mbox{ exists for every } \beta > \frac{1-\kappa}{\kappa}.
$$
On the other hand, by Theorem \ref{T:Fatou}, 
for $\nu$-a.e. $z \in \partial D$
$$
\lim_{ A^{\beta}_z \ni x \rightarrow z} \frac{f(x)}{h(x)} ~\mbox{ exists for every } \beta > \frac{1-\kappa}{\kappa}.
$$
This proves that for $\nu$-a.e. $z \in \partial D$,
$$
\lim_{ A^{\beta}_z \ni x \rightarrow z} \frac{u(x)}{h(x)}~=~\lim_{ A^{\beta}_z \ni x \rightarrow z} \frac{f(x)}{h(x)}\,+\,\lim_{ A^{\beta}_z \ni x \rightarrow z} \frac{g(x)}{h(x)}  ~\mbox{ exists for every } \beta > \frac{1-\kappa}{\kappa}.
$$
In particular,  for $\nu$-a.e. $z \in \partial D$,
$$
\lim_{ A^{\beta}_z \ni x \rightarrow z} \frac{k(x)}{h(x)} ~\mbox{ exists for every } \beta > \frac{1-\kappa}{\kappa}.
$$
Moreover, the above limit is strictly positive  by (\ref{KM}). Therefore,  for $\nu$-a.e. $z \in \partial D$, we have
$$
\lim_{ A^{\beta}_z \ni x \rightarrow z} \frac{u(x)}{k(x)}~=~\lim_{ A^{\beta}_z \ni x \rightarrow z} \frac{u(x)/h(x)}{k(x)/h(x)} ~\mbox{ exists for every } \beta
 > \frac{1-\kappa}{\kappa}.
$$
\qed

\medskip

\begin{remark}
{\rm If $D$ is a bounded Lipschitz open set, then for every nonnegative  $(\mu, F)$-harmonic function $u$,
\begin{equation}\label{FF_FK}
\lim_{ A^{\beta}_z \ni x \rightarrow z} ~\frac{u(x)}{\int_{\partial D} K_D(x, w) \sigma (dw)}~ \mbox{ exists and is finite for }\sigma\mbox{-a.e. } z \in \partial D.
\end{equation}}
\end{remark}

\medskip

Now suppose that $n=2$, $D=B:=B(0,1)$, $x_0=0$ and $\sigma_1$ is the normalized surface measure on $\partial B$. Recall that a curve $C_0$ is called a tangential curve in $B$ which ends on $\partial B$ if $C_0 \cap \partial B = \{w_0\} \in \partial B$, $C_0 \setminus \{w_0\} \subset B$ and there are no $r > 0$ and $ \beta >1$ such that $C_0 \cap B(w_0,r) \subset A^{\beta}_{w_0} \cap B(w_0,r)$. Because of (\ref{KM}), we can show that our result is the best possible one using the same argument as in Lemma \ref{counter_a} and Theorem \ref{counter_thm}.

\medskip
\begin{lemma}\label{counter_a_FK}
Let 
$$k(x):=\int_{ \partial B } K_B(x, w) \sigma_1(dw).$$
Suppose $U$ is a measurable function on $\partial B$ such that $0 \le U \le 1.$ Let 
$$u(x):=\int_{ \partial B } K_B(x, w) U(w) \sigma_1(dw) = \frac1{2 \pi} \int^{2 \pi}_0 K_B(x, e^{i \theta}) U( e^{i \theta}) d  \theta$$
where $x \in B.$
Suppose that  $0 < \lambda < \pi $ and $U(e^{i \theta})=1$ for $\theta_0 - \lambda \le \theta \le  \theta_0 + \lambda.$
Then there exists a $\delta = \delta (\eps, \alpha)$ such that 
$$
1-\eps \le \frac{u(\rho e^{i \theta_0})}{k(\rho e^{i \theta_0})} \le 1  ~~~~~\mbox{if } \rho > 1- \lambda \delta.
$$ 
\end{lemma}

\medskip

\begin{thm}\label{counter_thm_FK}
 Let $$k(x):=\int_{ \partial B } K_B(x, w) \sigma_1(dw).$$ Let $C_0$ be a tangential curve in $B$ which ends on $\partial B$ and let $C_{\theta}$ be the rotation of $C_{0}$ about $x_0$ through an angle $\theta.$
Then there exists a positive $(\mu, F)$-harmonic function $u$ in $B:=B(x_0,1)$ such that for a.e. $\theta \in [0, 2\pi]$ with respect to Lebesgue measure,
$$
\lim_{|x| \rightarrow 1, x \in C_{\theta}} \frac{u(x)}{k(x)} \mbox{ does not exist}.
$$
\end{thm}
\medskip

Now we assume  $D$ is a bounded $C^{1,1}$-open set and let 
$$
\psi (r):= 2^{-(n+\alpha)} \, \Gamma \left( \frac{n+\alpha}{2} \right)^{-1}\, 
 \int_0^\infty s^{\frac{n+\alpha}{ 2}-1} 
e^{-\frac{s}{ 4} -\frac{r^2}{ s} } \, ds,
$$
which is a smooth function of $r^2$, 
and $m>0$ be a constant.  We define
\begin{eqnarray*}
&&K^m_t := \exp \left(\sum_{0<s \le t} \ln(1+F_m(X^D_{s-}, X^D_s)) - {\cal A} (n, -\alpha)\int^t_0 \int_D F_m(X^D_s, y)|X^D_s-y|^{-\alpha-n} dy ds\right.\\ 
&&\left.~~~~~~~~~~~~~~~~- \int_0^t q(X^D_s) ds\right)
\end{eqnarray*}
where 
\begin{eqnarray*}
&&F_m(x,y)=\psi(m^{1/\alpha} |x-y|)-1, ~~
{\cal A} (n, \, -\alpha)=\frac{\alpha 2^{\alpha-1}\Gamma(\frac{\alpha+n}2)}
{\pi^{n/2}\Gamma(1-\frac{\alpha}2)}\\
 &\mbox{and}& 
q(x)={\cal A} (n, -\alpha)  \int_{D^c} F_m(x,y) dy.
\end{eqnarray*}

In Chen and Song \cite{CS4}, they obtained $X^m$, (killed) relativistic stable process in $D$ with parameter $m > 0 $ from $X^D$ through nonlocal Feynman-Kac transform $K^m_t$ for $\alpha \in (0,2)$. That is,
$$
 \E_x\left[f(X^m_t)\right]:= \E_x\left[f(X^D_t)K^m_t\right]
$$ 
for every positive Borel measurable function $f$, and 
$x \to \E_x[K^m_{\tau_D}]$ is bounded between
two positive constants.
Thus as a consequence of Theorem \ref{T:Fatou_FK}, we have the following.

\medskip

\begin{thm}\label{T:Fatou_R}
Let $D$ be a bounded $C^{1,1}$-open set and $\nu$ be a finite measure on $\partial D$. Define
$$
k(x) := \int_{\partial D} K_m(x, w) ~\nu (dw), ~~~~x\in D
$$ where  $K_m(x, w)$ is the Martin kernel for  $X^m$.
If $u$ is a nonnegative  harmonic function with respect $X^m$, then for $\nu$-a.e. $z \in \partial D$, 
$$
\lim_{ A^{\beta}_z \ni x \rightarrow z} \frac{u(x)}{h(x)} \mbox{ exists for every } \beta
 >1.
$$
\end{thm}

\medskip

\medskip

{\it Acknowledgment}: This paper is a part of the author's PhD thesis. He especially thanks Professor Z.-Q. Chen, his PhD thesis advisor, for his guidance and encouragement.
Thanks are also due to Professor R. Song for helpful comments. 
\vskip 0.5truein

\begin{singlespace}

\end{singlespace}

\bigskip
\bigskip

Panki Kim

Department of Mathematics,  
University of Illinois, 
Urbana, IL 61801, USA.

E-mails: pankikim@math.uiuc.edu

\end{doublespace}
\end{document}